\documentclass[a4paper,10pt]{article}

\newtheorem{thm}{Theorem}
\newtheorem{cor}[thm]{Corollary}
\newtheorem{rem}[thm]{Remark}
\newtheorem{lem}[thm]{Lemma}
\newtheorem{prop}[thm]{Proposition}
\def \a{{\alpha}}
\def \b{{\beta}}
\def \D{{\Delta}}
\def \DD{{\mathcal D}}

\def \e{{\varepsilon}}

\def \g{{\gamma}}
\def \G{{\Gamma}}

\def \l{{\lambda}}

\def \o{{\omega}}
\def \O{{\Omega}}
\def \p{{\varphi}}
\def \t{{\vartheta}}

\def \m{{\mu}}
\def \s{{\sigma}}

\def \A{{\cal A}}

\def \E{{\bf E}\, }

\def \P{{\bf P}}

\def \qq{{\qquad}}
\def \R{{\bf R}}

\def \T{{\bf T}}
\def \td{{\widetilde D}}
\def \ua{{\underline{a}}}

\def \ud{{\underline{d}}}

\def \uz{{\underline{z}}}

\def \Z{{\bf Z}}
\def \zz{{\cal Z}}

\def \noi{{\noindent}}

\def\qed{\hbox{\vrule height 6pt depth 0pt width
6pt}}
\def\cqfd{\hfill\penalty 500\kern 10pt\qed\medbreak}
\font\phh=cmcsc10
at  8 pt
\scrollmode
\title{Supremum of
 Random Dirichlet Polynomials
with Sub-multiplicative Coefficients}
 \author{
   Michel Weber}

\begin{document}

\maketitle

\begin{abstract}
We study the  supremum  of  random Dirichlet polynomials
$D_N(t)=\sum_{n=1}^N\e_n d(n) n^{- s }$, 
where $(\e_n)$ is a sequence of independent Rademacher random variables,
and  $ d  $ is a sub-multiplicative function.
The approach is gaussian and
  entirely based on   comparison properties of Gaussian processes, with no use of the metric entropy method.  
\end{abstract}


\section{Introduction}
   Let $\e= \{\e_n, n\ge 1\}$ denote  
a sequence of independent Rademacher random variables ($\P\{ \e_i=\pm 1\} =1/2$) defined on a basic
probability space $(\O, \A, \P)$.
Consider the random Dirichlet polynomials in which $s=\s+it $,
\begin{equation} \label{e12}
  \DD(s) =\sum_{n=1}^N \e_n d(n) n^{-s}.
\end{equation}  
     In  a recent work  \cite{LW2}, (see references therein for related results, notably Queffelec's works) we obtained sharp estimates of the supremum   of $
\DD(s)$,   under moderate growth condition on coefficients. Put
\begin{eqnarray} \label{char}
 D_1(M)&=& \sum_{m=1}^M d(m), \qquad\   \td_1(M)=\max_{1\le m\le M} \frac {D_1(m)}{m}   ,\cr 
  D_2(M) &=& \sum_{m=1}^M d(m)^2,  \qquad  \td_2^2(M)=\max_{1\le m\le M} \frac {D_2(m)}{m} .
\end{eqnarray}

We showed 
\begin{thm} \label{t1}
 Let $0\le\s \le 1/2$ and assume that
\begin{equation} \label{basic_d}
      d(k p^j) \le C d(k) j^H
\end{equation}
 for some positive $C,H$, any positive integer $k,j$ and any prime $p$.
 Then there exists a constant $C_{\s,d}$ depending on $d$ and $\s$ such that for any integer $N\ge 2$
 \begin{equation}\label{main_upper}
\E \sup_{t\in \R} \big|\DD(\s+it) \big| \le C_{\s,d}\, {  N^{1-\s}  \td_2(N)  \over \log N} \ .
\end{equation}
Moreover, if for some $b<1/ (\sqrt{5}+1)\approx 0.31$
\begin{equation}\label{extra}
\td_2(M) \le C  M^b,
\end{equation}
then
\begin{equation}\label {est.ex}
\E \sup_{t\in \R} \big| \DD(\s+it) \big| \le C_{\s,d}\, {  N^{1-\s}   \over \log N} \ .
\end{equation}
\end{thm}
  
Suppose $d(n)$ is a multiplicative function: $d(nm)\le d(n)d(m)$ if $n,m$ are
coprimes. Then   condition (\ref{basic_d}) is satisfied iff
\begin{equation} \label{basic_o} d(p^{r+j})\le C d(p^r) j^H,
\end{equation}
for some $C>0$, $H>0$ and any $j\ge 1$, $r\ge 0$. This last condition is fulfilled when  for
instance   
\begin{equation} \label{basic_m}
       \frac{d(p^{k+1})}{d(p^{k })} \le (1+\frac{1}{k})^H, \qquad  k=0,1,  \ldots
\end{equation}
a property which is    satisfied for a
relatively  wide class of multiplicative functions, among them,  the divisor function and 
 $ d_1(n)= \l^{\omega(n)}  $,
where  $\l>1$ and $\o(n)= \#\{p: p\mid n\}$ is    well-known additive prime divisor function.
\smallskip\par However,  condition (\ref{basic_d})  requires that  
     $ d(  p^j)= \mathcal O( j^H)$. Thus    Theorem  \ref{t1} does not apply to some classical  multiplicative  functions such as
$$  d_2(n)= \l^{\Omega(n)},$$
where   $\O(n)= \sum_{ p^\nu||n } \nu$  is the other prime divisor function.

  \bigskip\par 
The main concern of this  work is to   show that the approach used in \cite{LW2} can   be still adapted and further,   simplified,  
to obtain extensions  for a much larger class of multiplicative functions including these examples,   and also for sub-multiplicative
functions, namely functions satisfying the weaker condition:
\begin{equation}\label{hal0}
 d(nm)\le d(n)d(m) \qq \hbox{provided $(n,m)=1$.}
\end{equation} 
For instance,   $d(n)=e^{(\log n)^\a}$, $0<\a<1$ is   sub-multiplicative, as well as   function $d_K(n)= \chi\{(n,K)=1\}$ in Example
2. The related random Dirichlet polynomials are studied in this paper. 
\smallskip\par We obtain a
general upper bound, which also contains and improve the main results in
\cite{LW1},
\cite{LW2} (Theorem 1.1  and Theorem
\ref{t1} respectively). Introduce some notation. Let $2=p_1<p_2<\ldots$ be the sequence of all primes, and let $\pi(N)$ denote     the
number of prime numbers
  less or equal to $N$. The following decomposition   is   basic $$
\{2,\ldots,N\}=\sum_{j=1}^{ \pi(N)} E_j \qq {\rm where}\qq E_j=\big\{ 2\le n\le N :  P^+(n)=p_j\big\},
$$
  $P^+(n)$ being the largest prime divisor
of $n$. It is natural to disregard  cells $E_j$ such that $d(n)\equiv0$,  $n\in E_j$. We thus
set
$$\mathcal H_d=\big\{1\le j\le \pi(N):\  d_{| E_j}\not\equiv 0\big\}, \qquad  \tau_d=\max \mathcal (H_d).$$
 
 Consider now the following condition:   
 \begin{equation}\label{pest}p\mid n\  \Longrightarrow  \ d(n)  \le C\, d({n\over p }),   \qquad  {\rm and}\qquad    d(p^j)\le C_1\l^j,
 \end{equation}
for some   positive  $C,C_1, \l$ with $\l<\sqrt 2$,  any prime
 number $p$, any integers $n,j$.
  Clearly, if $C<\sqrt 2$, the second property is implied by the first. But this is not always so. Consider the following example. Fix
some prime number $P_1$ as well some reals $1<\l_1<\sqrt 2$, $C_1\ge 1$, and put 
\begin{equation} d  (n)= 
\cases{ C_1\l^j,  & \qq if $P_1^j||\, n$,\cr
 1,  & \qq if $(n,P_1)=1$.
}
\end{equation}
Then $d$ is sub-multiplicative, and satisfies condition (\ref{pest}) with a  constant $C$ which   has to be larger than $  C_1\l$. 

That $d$ be sub-multiplicative is easy: let $n,m$ be coprime integers. If $(n,P_1)=1$ and $(m,P_1)=1$,   then
$d(n)=d(m)=d(nm)=1$.  If $ P_1^j||n$ and $(m,P_1)=1$, then $d(n)=C_1\l^j$, $d(m)=1$; so $d(nm)=d(n)=d(n)d(m)$. Finally if $ P_1^j||n$
and $ P_1^k||m$, then $d(nm)= C_1\l^{\max(j,k)}\le C_1^2\l^{ j+k}=d(n) d(m)$.

Now let $p$ be such that $p|n$. If $p\not=P_1$, either $P_1\not| n$ and then $d(n)=d(n/p)=1$, or $P_1^j || n$ and  
$d(n)=d(n/p)=C_1\l^j$. If $p =P_1$, assume first  $P_1  || n$,
then 
$d(n/p)=1$, and in order that $d(n)=C_1\l \le Cd(n/p) $, one must take $ C\ge  C_1\l$.   
  Finally if $P_1^j  || n$ with $j\ge 2$, then $d(n)=C_1\l^j=\l d(n/p)\le Cd(n/p)$. 
 It remains to observe that $d(p^j)=1\le C_1\l^j$, if $p\not= P_1$; and by definition $d(P_1^j) = C_1\l^j$. This proves our claim.

More generally, let $P_1<\ldots <P_J$ be $J$ prime numbers, together with reals $C_1<\ldots <C_J$   and $\l_1<\ldots <\l_J$  
such that  $1<\l_j<\sqrt 2$ and $C_j\ge 1$ for all $j$, and form the corresponding   functions $d_1, \ldots d_J$. The product of
sub-multiplicative functions being again a sub-multiplicative function, we deduce that the product  $d=d_1. \ldots d_J$  is another
example of sub-multiplicative function satisfying condition (\ref{pest}), with a constant $C$ which has to be greater than $C_1\l_1\ldots
C_J\l_J$.
\medskip\par
 We prove
\begin{thm} \label{t2}
 Let   $d$ be a non-negative  sub-multiplicative   function. Assume that condition (\ref{pest}) is realized.
 Let $0\le\s < 1/2$. Then there exists a
constant $C_{\s,d}$ depending on $\s$ and $d$ only, such that   for any integer $N\ge 2$,
 $$
 \E\, \sup_{t \in \R} |\mathcal D (\s+it) |\le  C_{\s,d}  \, \td_2(N)\, B,$$
where $$
B=\cases{
  {N^{1/2-\s}   \tau_d^{1/2} \over (\log N)^{1/2}}
   &,\ {\rm if} \ $\big( {N\log\log N\over \log N}\big)^{1/2}  \le
\tau_d\le \pi(N) ,$  \cr
& \cr
  {N^{3/4-\s} (\log\log N)^{1/4}\over (\log N)^{3/4}}
   &,\ {\rm if} \  $\big( {N \over (\log N)\log\log N}\big)^{1/2}\le
\tau_d \le \big( {N\log\log N\over \log N}\big)^{1/2},$ \cr
& \cr
 N^{1/2-\s}        
          \big({ \tau_d \log\log\tau_d \over  \log\tau_d  }\big)^{1/2}
   &,\ {\rm if} \  $1\le \tau_d \le
\big( {N \over (\log N)\log\log N}\big)^{1/2} .$
} $$
 \end{thm}
 \medskip\par 
Observe that   condition  (\ref{basic_d}) implies condition  (\ref{pest}).    Indeed, write
$n=kp$ and take $j=1$. We get
$
       d(n)=d(k p ) \le C d(k)=   C d(n/p)
$. 
 Fix some real $\l$, $1<\l<\sqrt 2$. Then    $d(  p^j) \le C d(1) j^H \le C_1\l^j$, for some suitable constant $C_1$. 
\smallskip\par
Further, function $d_1$ obviously
satisfies condition (\ref{pest}), whereas we know that it does not satisfy condition (\ref{basic_d}). 
\smallskip\par
The bounds given in Theorem 
\ref{t2} being all less than $C_{\s,\l}\td_2(N)  \, {N^{1 -\sigma} \over  \log N  }$, we therefore deduce that  Theorem  \ref{t1} is
strictly included in Theorem 
\ref{t2}. We give two classes of examples of application.
\bigskip\par  
 {\it Example 1.} Consider   multiplicative functions satisfying the following   condition:
\begin{equation} \label{basic_n}
      \frac{d(p^{a})}{d(p^{a-1 })} \le \l, \qquad  a= 1,2, \ldots
\end{equation}
Clearly (\ref{basic_n}) is strictly weaker than
(\ref{basic_m}). Further it implies   (\ref{pest}). First $d(p^j)\le d(1)\l^j $. Next, let $p\mid   n$ and $a$ denote the $p$-valuation
of
$n$:  $p^a||n$. By multiplicativity of
$d(.)$ and condition (\ref{basic_n})
\begin{eqnarray*} d(n) =        d({n\over p ^{ a }}  )   \, d(   
p ^{ a  } )
= d({n\over p ^{
a }}  )\,d(     p^{ a -1 })\, {d(    p^{ a  })\over d(    p^{ a -1 } )}=
  d({n\over p })    {d(    p ^{ a  })\over d(    p ^{ a -1} )}\le \l \, d({n\over p }) .
\end{eqnarray*}
Thus (\ref{pest}) is fulfilled.
 Notice that (\ref{basic_n}) implies
\begin{equation}  \label{bounded}
M_d:=\sup_{p} d(p)   <\infty
\end{equation}
with $M_d\le \l d(1)$. 
\medskip\par
Under condition (\ref{basic_n}), estimates for $\td_2(N)$ are known.  By theorem 2 of
\cite{HR} (see also
\cite{Hi}), any non-negative multiplicative function $d$ satisfying a Wirsing type condition
\begin{equation}\label{hr_0}
 d(p^m)\le \l_1\l_2^m,
\end{equation}
for some constants $\l_1>0$ and $0<\l_2<2$ and all prime powers $p^m\le x$, also satisfies
\begin{equation}\label{hr_1}
 {1\over x}\sum_{n\le x} d(n)\le C(\l_1,\l_2) \exp\Big\{\sum_{p\le x} {d(p)-1\over p}\Big\},
\end{equation}
where $C(\l_1,\l_2)$ depends on $ \l_1,\l_2 $ only.
 \smallskip\par As $d$ satisfies   (\ref{basic_n}), if $\l<\sqrt 2$,   condition (\ref{hr_0}) is   verified with $\l_1=M_d$, $ \l_2=\l$. Since $d^2$ is
multiplicative and satisfies (\ref{basic_n}) with  $\l^2<2$, we also have  that $d^2$ verifies condition
(\ref{hr_0}) as well.   Consequently,     from (\ref{hr_1}) follows that   
\begin{eqnarray}\label{hr_2}
  \td_1 (N)  & \le & C(\l ) \exp\Big\{\sum_{p\le N} {d(p)-1\over p}\Big\}\cr
\td_2 (N)  & \le & C(\l ) \exp\Big\{\sum_{p\le N} {d^2(p)-1\over p}\Big\}, 
\end{eqnarray}
for some constant
$C(\l )$ depending on $ \l  $ only. Recall that there exists an absolute constant $c_1$ such that for $x\ge 2$
\begin{equation}\label{formula}
\Big|\sum_{p\le x} {1 \over p} -\log\log x -c_1\Big|<{5\over \log x} . \end{equation}
Thus 
 $$\sum_{p\le x} {d(p) \over p} \le M_d\sum_{p\le x} {1 \over p}   \le M_d \log (c_2\log x) $$
and similarly
 $$ \sum_{p\le x} {d^2(p) \over p} \le  M_d^2\log (c_2\log x).$$

Thereby under condition (\ref{basic_n}), we have the following estimates
\begin{equation}\label{hr_3}
  \td_1 (N)   \le  C(\l )(\log N)^{M_d}
,\qq \td_2 (N)    \le   C(\l ) (\log N)^{ M_d^2}. 
\end{equation}
  For   functions $d_1,d_2$, there is also a  standard
way to proceed. Letting $\tau =\pi(N)$, we have for $d_2$ for instance
$${1\over N} \sum_{n\le N} \l^{\Omega(n)}\le  \sum_{n\le N} {\l^{\Omega(n)}\over n} \le \sum_{\a_1=0}^\infty\ldots
\sum_{\a_\tau=0}^\infty{\l^{\a_1+\ldots+\a_\tau}\over p_1^{\a_1}\ldots p_\tau^{\a_\tau}}=\prod_{j=1}^\tau\Big(1-{\l\over
p_j}\Big)^{-1}$$
 which can be evaluated by means of (\ref{formula}).
\medskip\par The restriction $\l<2$ can be relaxed into $\l<  q$, when considering, instead
of  
$
\DD(s)$,   random Dirichlet polynomials  based on sets of integers having all their prime divisors greater or equal to $q$, e.g. on some 
arithmetic progressions. To go
beyond a condition of type (\ref{basic_n}), notably to work under the weaker  condition (\ref{hr_0}), one has   probably to perform another approach than the one
based on a decomposition into random processes as appearing in (\ref{dec}) below.

 \medskip\par
{\it Example 2.}  Take some positive integer
$K$, 
  and put
$$
 d_K (n)= 
\cases{ 1,  & \qq if $(n,K)=1$\cr
 0,  & \qq if $(n,K)>1$.
}
 $$
Then $d_K$ is   sub-multiplicative.  Let $p\, | n$. By definition, $d_K(n/p)=0$ iff
$(n/p,K)>1$, in which case $(n,K)>1$ and so $d_K(n)=0$. Thus $d_K(n)\le d_K(n/p)$. Now if $d_K(n/p)=1$, that $d_K(n)\le d_K(n/p)$ is
trivial. Besides $d_K(p^j)=d_K(p)\le 1$. Therefore condition (\ref{pest}) is satisfied  with $C=1=\l$. And by (\ref{e12}),  this   
defines the  remarkable class of random Dirichlet polynomials,
\begin{equation}
  \DD(s) =\sum_{  (n,K)=1\atop 1\le n\le N } \frac{\e_n }{n^{ s}}  ,
\end{equation}  
which  naturally extends   the one of ${\cal E}_\tau$-based Dirichlet
polynomials   considered in \cite{Q3} and \cite{LW1}. Indeed,  
recall that ${\mathcal E}_\tau=\big\{2\le n\le N  : P^+(n)\le p_\tau\}$. Define 
 
\begin{equation}\label{ktau} K_\tau=
\cases{   \prod_{\tau<\ell\le \pi(N)}
p_\ell   & \textrm{ if }$ \tau<\pi(N)$\cr
   1   & \textrm{ if }  $\tau=\pi(N)$ .} 
\end{equation} Then
 $n\in {\cal E}_\tau$, $n\le N$, iff $(n,K_\tau)=1$, namely $d_{K_\tau}(n)=1$. So that
\begin{equation}\label{relktau}\sum_{n\in {\cal E}_\tau } \frac{\e_n }{n^{ s}}=\sum_{  n=1}^N d_{K_\tau}(n)\frac{\e_n }{n^{ s}}. 
\end{equation}
Consequently, the ${\cal E}_\tau$-based Dirichlet
polynomials are one example of Dirichlet
polynomials with sub-multiplicative weights. Here    $\mathcal H_{d_{K_\tau}}=\sum_{j\le \tau} E_j $.  
We therefore neglect   cells $E_j$, $j>\tau $. Further, we have $\td_1(N) =\td_2(N) \le 1$. 
\smallskip\par If we now specify Theorem \ref{t2} to this case, we get  
\begin{cor} \label{cor1} Let $0<\s<1/2$. We have 
\begin{equation} \label{newtau}
 \E\, \sup_{t \in \R} \Big|\sum_{n\in {\cal E}_\tau } \frac{\e_n }{n^{ \s +it}} \Big|\le  C_{\s}   \, B,\qq {\rm where}
 \end{equation}
$$B=\cases{
  {N^{1/2-\s}   \tau ^{1/2} \over (\log N)^{1/2}}
   &,\ {\rm if} \ $\big( {N\log\log N\over \log N}\big)^{1/2}  \le
\tau \le \pi(N) ,$  \cr
& \cr
  {N^{3/4-\s} (\log\log N)^{1/4}\over (\log N)^{3/4}}
   &,\ {\rm if} \  $\big( {N \over (\log N)\log\log N}\big)^{1/2}\le
\tau  \le \big( {N\log\log N\over \log N}\big)^{1/2},$ \cr
& \cr
 N^{1/2-\s}        
          \big({ \tau \log\log\tau  \over  \log\tau  }\big)^{1/2}
   &,\ {\rm if} \  $1\le \tau  \le
\big( {N \over (\log N)\log\log N}\big)^{1/2} .$
} $$
\end{cor}

By comparing this with the upper bound part
of Theorem 1.1 in
\cite{LW1},
 we observe that the bounds obtained are either identical (if $N^{1/2}\le \tau\le \pi(N)$), or strictly better. For
instance,    when $ ( {N \over (\log N)\log\log N} )^{1/2}\le
\tau\le
\big( {N\log\log N\over \log N}\big)^{1/2}$, we have 
$$ {N^{3/4-\s} (\log\log N)^{1/4}\over (\log N)^{3/4}}\ll {N^{3/4-\sigma}
\over (\log N)^{1/2}},$$
thereby yielding a better bound.  
\medskip\par 
When the order of $\tau $ is small, we will prove the following strenghtening in which    $N$ disappears from the estimates. Put 
$$\Pi_\s(\tau)= \prod_{\ell =1}^\tau \big[{1\over 1-   p_\ell^{-2\s }} \big]. $$

 \begin{thm} \label{smalltau} Assume that $\tau=o(\log N)$. Let $0<\s<1/2$. Then, there are    $c_\s,C_\s$ depending on $\s$ only, such that
\begin{equation}\label{estsmalltau}
c_\s \, {\Pi_\s(\tau)^{1/2}\,\tau^{1-\s}  \over( \log \tau)^\s }\le \E \sup_{t\in \R}\big|\sum_{  n\le N\atop P^+(n)\le p_{\tau}}  {\e_n
\over n^{ \s+ it}}\big|
 \le C_\s \,    \Big({\Pi_\s(\tau)^{1/2}\, \tau^{{3\over 2}- 2\s} \over(\log \tau)^{2 \s} }
    \Big).
\end{equation} 
And if $\s=1/2$, there are absolute constants $C_1,C_2$ such that
$$    C_1    {\tau    }^{1/2} \le \E \sup_{t\in \R}\big|\sum_{j=1}^\tau \sum_{n\in E_j}  {\e_n\over \sqrt n}n^{-it}\big|\le   C_2   
{\tau    }^{1/2} (\log\log  \tau )^{1/2} .$$
\end{thm}
 \bigskip\par
  Let now $K$ be unspecified. There is no loss to assume $K$ is squarefree. First examine  the case when $K$ has few prime
divisors. Suppose
 \begin{equation}\label {exmp}
  \sum_{  p|K\atop   p\le N }  
  \frac{1 }{p^{ \s}} =o({  N^{1-\s}   \over \log N}). \end{equation}   
  Using Bohr's lower bound
  \begin{equation}\label {exm2}
 \E \sup_{t\in \R} \Big| \sum_{  (n,K)=1\atop 1\le n\le N } \frac{\e_n }{n^{ s}}  \Big| \ge C   \sum_{  (p,K)=1\atop   p\le N }
  \frac{1 }{p^{ \s}}     .
 \end{equation}
 
We get with \ref{t2}  a two-sided estimate 
 \begin{equation}\label {ex3}
 C\ {  N^{1-\s}   \over \log N} \le \E \sup_{t\in \R} \Big| \sum_{  (n,K)=1\atop 1\le n\le N } \frac{\e_n }{n^{ s}}  \Big| \le C\ { 
 N^{1-\s}   \over \log N} \,  .
 \end{equation}

The case of a number $K$ with many prime divisors is more complicated. By the comment previously made, this concerns the case 
 \begin{equation}\label {exmpmany}
  \sum_{  p|K\atop   p\le N }  
  \frac{1 }{p^{ \s}} \asymp {  N^{1-\s}   \over \log N} . \end{equation}   
 We   restrict ourselve to integers $K$ of     type
$$K=\prod_{p|K\atop p\le p_\nu} p\cdot\prod_{   p_\nu<p\le N} p , $$
where $1\le \nu <\pi(N).$ This amounts to consider the random Dirichlet polynomials
$$\sum_{     1\le n\le N\atop (n,K)=1} \frac{\e_n }{n^{ s}}=\sum_{       n\in F_\nu \atop (n,K)=1} \frac{\e_n }{n^{ s}}.$$
We will assume $\nu$ to be  not too large. More precisely, we assume, in accordance with  Corollary \ref{cor1} 
$$  \nu \le
\big( {N \over (\log N)\log\log N}\big)^{1/2}.$$

\begin{thm} \label{t3} Let $0<\s<1/2$.   There exists a constant $C_\s$ depending on $\s$ only such that 
$$
\E \sup_{t\in \R} \Big| \sum_{   
(n,K)=1}
\frac{\e_n }{n^{ s}}\Big|
\le  C_\s  \, N^{1/2- \s}       
   \max\Big(1,\sum_{k\le \nu\atop
 p_k\not|K  } {1\over   {p_k} }   
    \Big)^{1/2}\Big[\sum_{k\le \nu\atop
 p_k\not|K  } {1\over \sqrt {p_k} }\Big] .
$$  
  \end{thm} 
 
   
\medskip\par
{\it Example 3.} Fix some integer $N\ge 1$, and consider the truncated divisor function
$$d_N(n)=\#\{ k\le N: k|n\}. $$ 
This function, which occurs    in many  important arithmetical questions,   is sub-multiplicative. Take $n$ and $m$ coprimes. If $k\le
N$ is such that
$k|mn
$, then $k$ is uniquely written   $k=k_1k_2$, $(k_1,k_2)=1$, $k_1|m$, $k_2|n$; and naturally $k_1\le N$, $k_2\le N$. We infer that
$d_N(mn)\le d_N(m)d_N(n)$. 
\smallskip\par Further, it satisfies our condition (\ref{pest}).
Let $p|n$, if $p>N$ then $d_N(n)=d_N({n\over p})$. Otherwise, if $p\le N$, let ${\cal K}=\{k\le N:(k,p)=1\}$.   For $k\in {\cal K}$
such that $k|n$,   the $p$-height $p(k)$ of $k$ denotes the largest integer $a$ so that
$p^ak|n$ and $p^ak\le N$. The divisors of $n$ are of
type
$p^ak$, $k\in {\cal K}$. Further if $ p^ak_1  =p^bk_2$, $k_1,k_2\in {\cal K}$, necessarily $k_1=k_2$. Indeed,   it is
obvious if $a=b$; and if
$a>b$  we get $p|k_2$, which excluded.  Consequently 
$$d_N(n)= \sum_{k\in {\cal K}\atop k|n} (1+p(k)), \qq d_N({n\over p})= \sum_{k\in {\cal K}\atop k|n} [1+(p(k)-1)^+]. $$
As   for any integer $a\ge 0$, $ 1+a\le 2[1+(a-1)^+]$, we deduce
$$d_N(n)\le 2 d_N({n\over p}). $$
  And choosing any $\l>1$, we obviously have  
 $d_N(p^j) =\#\{ \ell\le j: p^\ell\le N\}\le j\le C\l^j
 $.


\section{Proof of Theorem \ref{t2}.}
Although the proofs are much in the spirit of   proofs of the main results in \cite{LW1},\cite{LW2}, 
there are substancial changes and simplifications. First, we work
from the beginning with Gaussian processes. Further, the delicate step of estimating  some related metric and computing associated entropy numbers is   notably
simplified.  Cauchy-Schwarz's inequality and the comparison properties of Gaussian processes indeed   allow to avoid any computation (see before
(\ref{e211})), and also give rise to   strictly better estimates. 

This   further allowed us to consider random Dirichlet polynomials with
more complicated arithmetical structure, like the one of "Hall type" built from the  sub-multiplicative functions $d_K$, where entropy
numbers seem hard to estimate efficiently.

 \medskip\par  Let $\tau=\tau_d$ and let $a_j(n)$ denote  the $p_j$-valuation of integer $n$.     Put  
$$\ua(n)=\big\{a_j(n), 1\le j\le \tau\big\},\qq\qq (n\le N).$$  
   Let also $\T=[0,1[=\R/\Z$ be the torus.
A first  classical reduction  allows to replace the Dirichlet polynomial by some relevant   trigonometric
polynomial.  To any Dirichlet polynomial $P(s)=\sum_{n=1}^N a_n n^{-s}$, associate the trigonometric
polynomial $ Q(\uz)$ defined   by
$$
 Q(\uz)= \sum_{n=1}^N  a_n n^{-\s}e^{2i\pi\langle \ua(n),\uz\rangle}, \qq \uz= (z_1,\ldots, z_\tau) \in \T^\tau.
$$
    
By Kronecker's Theorem (\cite{HW}, Theorem 442)  
\begin{equation} \label{e11}
 \sup_{t\in \R} \big|P(\s +it)\big| =\sup_{\uz \in \T^\tau} \big|Q(\uz)\big|,
\end{equation}
 as observed in \cite{B}. 
 \smallskip\par

\begin{rem} \rm  Naturally,  no similar  reduction   occurs  when     considering    the supremum over a given bounded
interval $I$. 
  However, when the length of $I$ is of
 exponential size with respect to the degree of $P$, precisely when 
 $$|I|\ge  e^{    (1+\e) \o  N(\log  N\o) 
  \log N  },$$  
  the related supremum   becomes
comparable, for $\o$ large,  to the one taken on the real line, with an error term of order $\mathcal O(\o^{-1}) $. This is in turn a
rather general phenomenon due to   existence of "localized" versions of Kronecker's theorem;    and in the present case to   Tur\'an's estimate (see  
\cite{W3} for   a   slighly improved form of it, and
      references therein).
When the length is of sub-exponential order, the study still seems to belong to the field of
application of the general theory of regularity of stochastic processes. 
\end{rem}  

In the technical lemma below, we collected some  useful estimates,     which
already appeared in \cite{LW2},  and are easily deduced from  
  the fact that if $a_n$ are complex numbers and
$b\in
\mathcal C^1([1,x])$, then
\begin{equation}\label{simple}\sum_{1\le n\le x} a_nb(n)= A(x)b(x)-\int_1^x A(t) b'(t) dt,
\end{equation}
where we let $A(t)=\sum_{n\le t}a_n$.   \begin{lem} \label{lem1} Let
$M\le N$ and
$0<\s< 1/2$. Then  
\begin{equation} \label{abel0}
 \sum_{m\le M} {d(m)^2\over  m^{ 2\s}} \le C \td_2^2(M) M^{1-2\s}.
\end{equation} 
\begin{equation} \label{abel1}
 \sum_{m\le M} ({N\over m})^{1/2} \ (\log({N\over m}))^{-1/2} d(m)
 \le  C \td_1(M)  (N M)^{1/2} \ (\log({N\over M}))^{-1/2}.
\end{equation}  
 \begin{equation}  \label{abel2}
 \sum_{k\le M} {d(k)^2 \over k^{
2\s}  }\le C \td_2(M)^2 (M)^{1-2\s}.
\end{equation}
\end{lem}

 \bigskip\par
  Now we can pass to the proof of Theorem  \ref{t2}.  Fix some integer $\nu$ in $[1,\tau]$. We denote  
$$F_\nu  
 = \sum_{1\le j\le \nu} E_j, \qq  F^\nu=\sum_{  \nu<j\le \tau }E_j . $$ Consider as in \cite{LW1},\cite{LW2}    the
decomposition  
$Q=Q^\e_1+Q^\e_2$, where
\begin{eqnarray*}
Q^\e_1(\uz) &=& \sum_{  n\in F_{\nu }} \e_n d(n) n^{-\s}
e^{2i\pi\langle\ua(n),\uz\rangle}, \
\cr
Q^\e_2(\uz) &=& \sum_{  n\in F^{\nu }}
\e_n  d(n) n^{-\s}e^{2i\pi\langle \ua(n),\uz\rangle}.
\end{eqnarray*}
 By the contraction principle (\cite{K} p.16-17)
 \begin{equation}\label{Q2A}
\E\, \sup_{\uz \in \T^\tau}\big|Q^\e_i(\uz)\big|
\le  4 \ \sqrt{ \pi \over 2  }\
\E\, \sup_{\uz \in \T^\tau}\big|Q_i(\uz)\big|,\qq (i=1,2)
\end{equation}
where $ Q_i $ is the same process as $Q_i^\e $
except that the Rademacher random variables $\e_n$ are replaced by
independent
${\cal N}(0,1)$ random variables $\m_n$. Consequently, both the supremums of $Q_1$ and of $Q_2$ can be estimated, via
  their associated $L^2$-metric. \medskip\par 
      Assume first  $0<\s<1/2$. We will establish the two following estimates:
 \begin{equation}\label{Q1B}
 \E\, \sup_{\uz \in \T^{\tau }}\big|Q_1(\uz)\big|\le
 C
N^{1/2-\s} \td_2(N) \big({\nu \log\log \nu\over \log \nu} \big)^{1/2},
\end{equation}
and \begin{equation}\label{Q2D}
\E\, \sup_{\uz \in \T^\tau}\big|Q_2(\uz)\big|
\le      C\bigg(   N^{1/2-\s}  \td_2(N/p_\nu) \ {\tau^{1/2} \over (\log\tau)^{1/2}}+{   N^{1-\s} \td_1(N/p_\nu)\over \nu^{1/2} \log\nu}
\bigg) .
\end{equation} 

 \bigskip\par 
First, evaluate the supremum of
$Q_2$.   Writing  
\begin{eqnarray*}
Q_2(\uz)&=& \sum_{ \nu <j\le \tau}
e^{2i\pi z_j  }\sum_{n\in E_j}
\mu_n d(n) n^{-\s} e^{2i\pi\{\sum_{k\not = j} a_k(n)z_k+ [a_j(n)-1]z_j\}}
\cr
&=& \sum_{ \nu <j\le \tau} e^{2i\pi z_j  }\sum_{n\in E_j}
\mu_n  d(n)  n^{-\s}e^{2i\pi\left\{\sum_{k} a_k({n\over p_j}) z_k\right\}}
\end{eqnarray*}
 next developing, gives $$=\sum_{ \nu <j\le \tau}  \cos  2 \pi   z_j   \sum_{n\in E_j}
 \mu_n  {d(n) \over  n^{ \s}}\ \cos  2 \pi \sum_{k} a_k({n\over p_j}) z_k $$
 $$+i\sum_{ \nu <j\le \tau}   \sin  2 \pi   z_j   \sum_{n\in E_j}
 \mu_n  {d(n) \over  n^{ \s}} \cos  2 \pi \sum_{k} a_k({n\over p_j}) z_k   $$
 $$+i\sum_{ \nu <j\le \tau}  \cos  2 \pi   z_j  \sum_{n\in E_j}
 \mu_n  {d(n) \over  n^{ \s}} \sin  2 \pi \sum_{k} a_k({n\over p_j}) z_k $$
 $$-\sum_{ \nu <j\le \tau}    \sin  2 \pi   z_j   \sum_{n\in E_j}
 \mu_n  {d(n) \over  n^{ \s}} \sin  2 \pi \sum_{k} a_k({n\over p_j}) z_k $$
 with $ {n/ p_j}\le N/p_j<N/p_\nu\le N/2$. 
  Each piece is, up to a factor $1,i, -1$, one of the possible realizations of the
random process $X $ defined by 
\begin{equation}\label{dec}
X (\g) =\sum_{ \nu <j\le \tau} \a_j
\sum_{n\in E_j}  \mu_n  {d(n)\over  n^{ \s}} \b_{{n\over p_j}},
\qq \g\in \G,
\end{equation}
where
$\g =\big((\a_j)_{\nu <j\le \tau}, (\b_m)_{1\le m\le N/2}\big)
$
and
$$\G=\big\{ \g  :  |\a_j|\vee |\b_m|\le \!\!1,  \nu < \!j\le \tau,
1\! \le m\le N/2\big\}.
$$
Here indeed
 $$
\a_j= \a_j(\uz)=
\cases{  \cos (2\pi  z_j), &\cr
 \hbox{\rm or}  &\cr
 \sin (2\pi  z_j),&
}
\qquad \nu <j\le \tau;
$$
and 
$$
\b_m = \b_m(\uz) =
\cases{ \cos\left( 2\pi \sum_{k} a_k({m}) z_k \right),  &\cr
\hbox{\rm or}  &\cr
\sin\left( 2\pi \sum_{k} a_k({m}) z_k \right),&
}
\qquad 1\le m\le {N\over 2} . 
$$
  
Consequently
\begin{equation}\label{Q2A}\sup_{\uz \in \T^\tau}\big|Q_2(\uz)\big|
\le
4\sup_{\g \in\G}  \big|X (\g)\big|.
\end{equation}
  The problem now reduces to estimating the
supremum over $\Gamma$ of the real valued Gaussian process $X$. We observe that
\begin{eqnarray*}
\|X_\g-X_{\g'}\|_2^2
&=&
\sum_{ \nu <j\le \tau} \sum_{n\in E_j}
 d(n)^2 n^{-2\s} \big[\a_j\b_{{n\over p_j}}-\a'_j \b'_{{n\over p_j}}\big]^2
\cr
&\le&
2\!\sum_{ \nu <j\le \tau}\sum_{n\in E_j}  d(n)^2 n^{-2\s} \big[(\a_j -\a'_j)^2+
(\b_{{n\over p_j}}-\b'_{{n\over p_j}})^2\big]  .
\end{eqnarray*}

 As $   p_j \mid n$, by   condition (\ref{pest}),
$d(n) \le \l \ d({n\over p_j})$; and so
\begin{eqnarray}
\sum_{ \nu <j\le \tau} \sum_{n\in E_j}  {d(n)^2 \over n^{ 2\s}} (\a_j -\a'_j)^2
&\le&
 \l^2 \sum_{\nu<j\le\tau} (\a_j -\a'_j)^2p_j^{-2\s}
\sum_{m\le N/p_j}  \frac{d(m)^2}{ m^{2\s}}
\cr
\cr  \label{e21}
&\le &
 \l^2  \sum_{\nu<j\le\tau} (\a_j -\a'_j)^2 {N^{1-2\s}\td _2^2(N/p_j) \over p_j},
\end{eqnarray}
where we used estimate (\ref{abel0}) of Lemma \ref{lem1}.

 Besides, by condition (\ref{pest}) again, we obtain
\begin{eqnarray}
\sum_{\nu <j\le \tau} \sum_{n\in E_j}
{   d(n)^2 (\b_{{n\over p_j}}- \b'_{{n\over p_j}})^2\over n^{2\s}}
&\le&
C\l^2 \sum_{m \le N/p_\nu}
(\b_m-\b'_m)^2 \big(\sum_{{\nu<j\le \tau\atop  mp_j\le N}}
{ d(m)^2 \over (mp_j)^{2\s}} \big)
\cr                                     \label{e22}
&:=&
C \l^2\sum_{m \le N/p_\nu} K_m^2
(\b_m-\b'_m)^2.
\end{eqnarray} 
   Let  $k\in (\nu,\tau]$ be such that
$N/p_k< m \le N/p_{k-1}$. Since $p_j \sim j   \log j$,
 we have
\begin{eqnarray*}
K_m^2 &=& \sum_{\nu<j\le k-1 }  d(m)^{2} (mp_j)^{-2\s}
\le
d(m)^{2} m^{-2\s}\ \sum_{j\le k-1 } p_j^{-2\s}
\cr
&\le&
C\ d(m)^{2} m^{-2\s}\ \sum_{j\le k } (j\log j)^{-2\s}
\le
C\  d(m)^{2} m^{-2\s}\  {k^{1-2\s} \over (\log k)^{2\s}}
\cr
&\le&  C d(m)^{2} m^{-2\s} \  {k \over p_k^{2\s}}
 \le
C m^{-2\s}  d(m)^{2} \  {k \over (N/m)^{2\s}}
\cr
&=& C \  d(m)^{2} {k \over N^{2\s}}\ .
\end{eqnarray*}
We have $k\log k\le Cp_k\le C\ {N\over m}$, and so
 $ k\le C\ {N\over m} \ (\log({N\over m}))^{-1}$.
 Thus 
\begin{equation} \label{e23a}
K_m\le C\  d(m)   N^{-\s}   ({N\over m})^{1/2} \ (\log({N\over m}))^{-1/2}
\ .
\end{equation}
By using estimate (\ref{abel1}) of Lemma \ref{lem1}
 \begin{eqnarray}
\sum_{m\le N/p_\nu} K_m
&\le& C\  N^{-\s}\ \sum_{m\le N/p_\nu}
({N\over m})^{1/2} \ (\log({N\over m}))^{-1/2} d(m)
\cr \label{e23b}
&\le&
{C N^{1-\s} \td_1(N/p_\nu)\over \nu^{1/2} \log\nu}\ .
\end{eqnarray}
Now define a second Gaussian process by putting for all $\g\in \G$
$$
Y(\g) =
\sum_{\nu<j\le\tau} \big({\td_2^2(N/p_j) N^{1-2\s}\over p_j}\big)^{{1/ 2}}\a_j\xi'_j
+
\sum_{m\le N/p_\nu} K_m \ \b_m\xi''_m
 :=
 \ Y'_\g + Y''_\g ,
$$
where $\xi'_i  $, $\xi''_j$  are independent ${\cal N}(0,1)$ random
variables. It follows from (\ref{e21}) and  (\ref{e22}) that for some suitable constant
$C$, one has the comparison relations: for all $\g, \g'\in \G$,
$$
\|X_\g-X_{\g'}\|_2\le C \|Y_\g-Y_{\g'}\|_2.
$$

By  the Slepian comparison lemma (\cite{L}, Theorem 4 p.190), since
$X_0=Y_0=0$,
we have
\begin{equation}\label{Q2B}
\E\, \sup_{\g\in \G}     |X_\g|\le
2 \E\, \sup_{\g\in \G}      X_\g \le
2 C \E\, \sup_{\g\in \G} Y_\g \le
2 C \E\, \sup_{\g\in \G} |Y_\g|.
\end{equation}
And with (\ref{Q2A}) 
\begin{equation}\label{Q2C}
\E\, \sup_{\uz \in \T^\tau}\big|Q_2(\uz)\big|
\le   C  
\E\, \sup_{\g \in \G}\big|Y(\g)\big|.
\end{equation}
It remains to evaluate the supremum  of $Y$. First of all,
$$
\E\, \sup_{\g\in \G}|Y'(\g)|
\le N^{{1\over 2}- \s}\sum_{\nu<j\le\tau}   p_j^{-1/2} \td_2(N/p_j).
$$
As $p_j\sim j\log j$, we have
$$
\sum_{\nu<j\le\tau}   p_j^{-1/2}
\le
   \sum_{1<j\le\tau}   p_j^{-1/2}
\le
 {C \tau^{1/2} \over (\log\tau)^{1/2}}\ ,
$$
thus
\begin{equation} \label{e24}
 \E\, \sup_{\g\in \G}|Y'(\g)|
 \le C\ N^{{1\over 2}-\s}  \td_2(N/p_\nu) \ {\tau^{1/2} \over (\log\tau)^{1/2}}
  \ .
\end{equation}

\noi To control the supremum of $Y''$, we use our
estimates for the sums of $K_m$ and write that
 \begin{equation} \label{e25}
 \E\, \sup_{\g\in \G} |Y''(\g)|
\le \sum_{m\le  N/p_\nu} K_m
 \le
{C  N^{1-\s} \td_1(N/p_\nu)\over \nu^{1/2} \log\nu} \ .
\end{equation}
 Therefore by reporting (\ref{e24}), (\ref{e25}) into (\ref{Q2C}), we get (\ref{Q2D}).

 \bigskip

Now, we turn to the supremum of $Q_1 $. Introduce the auxiliary Gaussian process
$$
\Upsilon (\uz) =
\sum_{n\in F_\nu } d(n)  n^{-\s}
\big\{\t_n \cos 2\pi \langle \ua(n),\uz\rangle +
\t_n'\sin 2\pi \langle \ua(n),\uz\rangle \big\}
,\qq  \uz\in \T^{\nu },
$$
where $\t_i$, $\t'_j$  are independent ${\cal N}(0,1)$
random variables. By symmetrization (see e.g. Lemma 2.3 p. 269 in \cite{PSW}),
\begin{equation}\label{Q1A}
\displaystyle{\E\, \sup_{\uz \in \T^{\nu }}\big|Q_1(\uz)\big|\le
\sqrt{8\pi}
\E\, \sup_{\uz \in  \T^{\nu }}\big|\Upsilon (\uz)\big|}.
\end{equation}
 
Further \begin{eqnarray*}
 \,\|\Upsilon(\uz)-\Upsilon(\uz)\big\|_2^2
&=& 4 \sum_{n\in F_\nu } {d(n)^2\over n^{2\s}} \sin^2(\pi \langle\ua(n),\uz -\uz'\rangle)
\cr
&\le &
4\pi^2\ \sum_{n\in F_\nu }  {d(n)^2\over n^{2\s}} |\langle\ua(n),\uz -\uz'\rangle|^2
\cr
&\le & 4\pi^2\ \sum_{n\in F_\nu }  {d(n)^2\over n^{ 2\s}}
    \Big[ \sum_{j=1}^\nu a_j(n) |z_j - z'_j|\Big]^2
 .
\end{eqnarray*}
Now 
$$\sum_{n\in F_\nu }  {d(n)^2\over n^{ 2\s}}\Big[ \sum_{j=1}^\nu a_j(n) |z_j - z'_j|\Big]^2=\sum_{j =1}^\nu  
      |z_{j } - z'_{j }|^2\sum_{n\in F_\nu }{a_{j }(n)^2d(n)^2\over n^{ 2\s}}$$$$ +\sum_{1\le j_1,j_2\le \nu\atop j_1\not=j_2}    |z_{j_1} - z'_{j_1}|\
|z_{j_2} - z'_{j_2}|\sum_{n\in F_\nu }{a_{j_1}(n)
    a_{j_2}(n)d(n)^2\over n^{ 2\s}}
 :=   S+R.$$

Examine first the contribution of the   rectangle  terms. Only those integers $n$ such that $ a_{j_1}(n) \ge 1$
   and $  a_{j_2}(n)\ge  1$ are to be considered.   Using sub-multiplicativity, we have
 \begin{eqnarray}\label{debutR} R  &\le&  \ \sum_{1\le j_1,j_2\le \nu\atop j_1\not=j_2}  |z_{j_1} - z'_{j_1}|\ |z_{j_2} -z'_{j_2}|
    \sum_{b_1,b_2=1}^\infty b_1 b_2 \sum_{{ n\le N, a_{j_1}(n)=b_1,\atop a_{j_2}(n)=b_2}} {d(n)^2\over n^{ 2\s}}
 \cr&\le &  C\sum_{1\le j_1,j_2\le \nu\atop j_1\not=j_2} |z_{j_1} - z'_{j_1}|  |z_{j_2} - z'_{j_2}|
    \sum_{b_1,b_2=1}^\infty      { b_1d(p_{j_1}^{ b_1})^2 \over  p_{j_1}^{2 b_1\s}   } {b_2 d( p_{j_2}^{ b_2})^2\over  p_{j_2}^{2 b_2\s}
}
 \cr& &  \qquad\times
\Big[\sum_{k\le {N \over p_{j_1}^{ b_1} p_{j_2}^{ b_2}}} {d(k)^2 \over k^{
2\s}  }\Big].
\end{eqnarray}
Examine now the contribution of the   square terms. We have
\begin{eqnarray}\label{debutS}S&\le& \sum_{j =1}^\nu  
      |z_{j } - z'_{j }|^2\sum_{b=1}^\infty \sum_{n\in F_\nu\atop a_{j }(n)=b  }{b^2d(n)^2\over n^{ 2\s}} \cr& \le &\sum_{j =1}^\nu  
      |z_{j } - z'_{j }|^2\sum_{b=1}^\infty {b^2d(p_j^b)^2\over p_j^{ 2b\s}}\sum_{m\le {N\over p_j^b}   }{ d(m)^2\over m^{ 2\s}} .
 \end{eqnarray}
 By estimate (\ref{abel2}) of Lemma \ref{lem1},  we have
 \begin{equation}  
 \sum_{k\le  {N \over p_{j_1}^{ b_1} p_{j_2}^{ b_2}}} {d(k)^2 \over k^{
2\s}  }\le C \td_2(N)^2 \Big[{N \over p_{j_1}^{ b_1} p_{j_2}^{ b_2}}\Big]^{1-2\s}.
\end{equation}
Hence
$$R\le   C \td_2(N)^2 \! \sum_{1\le j_1,j_2\le \nu\atop j_1\not=j_2}\!\! |z_{j_1} - z'_{j_1}|  |z_{j_2} - z'_{j_2}|
    \sum_{b_1,b_2=1}^\infty      { b_1d(p_{j_1}^{ b_1})^2 \over  p_{j_1}^{2 b_1\s}   } {b_2 d( p_{j_2}^{ b_2})^2\over  p_{j_2}^{2 b_2\s}
}
\Big[{N \over p_{j_1}^{ b_1} p_{j_2}^{ b_2}}\Big]^{1-2\s} 
$$
$$=   C \td_2(N)^2 N^{1-2\s} \sum_{1\le j_1,j_2\le \nu\atop j_1\not=j_2} |z_{j_1} - z'_{j_1}|  |z_{j_2} - z'_{j_2}|
    \sum_{b_1,b_2=1}^\infty      { b_1d(p_{j_1}^{ b_1})^2 \over  p_{j_1}^{ b_1 }   } {b_2 d( p_{j_2}^{ b_2})^2\over  p_{j_2}^{ b_2 }
}
. $$
But, by condition (\ref{pest}) 
 \begin{eqnarray} \label{pestelem} \sum_{b \ge 1}b{  d(p_j^{b})^2 \over   p_j^{b
 }  }&\le &  
 C  \sum_{b \ge 1}b\left({  \l^2   \over   2}\right)^{ b  } \left({  2   \over   p_j }\right)^{ b  }\le C  \left({  2   \over   p_j
}\right)\sum_{b \ge 1}b\left({ 
\l^2   \over   2}\right)^{ b  } 
\cr &\le &    { C_\l   \over   p_j }  
 .  
\end{eqnarray}
 From this follows that
\begin{equation}R\le    C_\l \td_2(N)^2 N^{1-2\s} \Big[\sum_{j =1}^\nu {|z_{j } - z'_{j }| \over p_j}
 \Big]^2
. \end{equation}

Further\begin{eqnarray}S&\le& C \td_2(N)^2\sum_{j
=1}^\nu  
      |z_{j } - z'_{j }|^2\sum_{b=1}^\infty {b^2d(p_j^b)^2\over p_j^{ 2b\s}}\Big[{N\over p_j^b}    \Big]^{1-2\s} \cr &\le& C 
\td_2(N)^2 N^{1-2\s}\sum_{j
=1}^\nu  
      |z_{j } - z'_{j }|^2\sum_{b=1}^\infty {b^2d(p_j^b)^2\over p_j^{ b }}\cr&\le &C \td_2(N)^2 N^{1-2\s}\Big[\sum_{j =1}^\nu  
     { |z_{j } - z'_{j }|^2\over p_j} \Big],
\end{eqnarray}
by arguing as in (\ref{pestelem}) for getting  the last inequality.
\smallskip\par
Consequently 
  \begin{equation} \label{e27}
 \|\Upsilon(\uz)-\Upsilon(\uz)\big\|_2\le C_\l N^{1/2-\s} \td_2(N)  \max \bigg(  \sum_{j=1}^\nu  {|z_{j} - z'_{j}|\over p_j}
, \Big[\sum_{j =1}^\nu  
     { |z_{j } - z'_{j }|^2\over p_j} \Big]^{1/2}\bigg)     .
\end{equation} 
\medskip

   We shall control the  Gaussian process
 $
\Upsilon  $ in a more simple  and more efficient way than  in \cite{LW1},\cite{LW2}. 
  By the Cauchy-Schwarz inequality 
\begin{eqnarray} 
   \sum_{j=1}^\nu  {|z_{j} - z'_{j}|\over p_j} &\le &    \Big( \sum_{j=1}^\nu    {|z_{j} - z'_{j}|^2\over p_j }\Big)^{1/2}
\Big(\sum_{j=1}^\nu  {1\over p_j }\Big)^{1/2}
\cr 
&\le &   \Big( \sum_{j=1}^\nu    {|z_{j} - z'_{j}|^2\over p_j }\Big)^{1/2} \Big(\sum_{j=1}^\nu  {1\over  j\log j }\Big)^{1/2}
\cr 
&\le &   (\log\log \nu)^{1/2}  \Big( \sum_{j=1}^\nu    {|z_{j} - z'_{j}|^2\over p_j }\Big)^{1/2} 
    .
\end{eqnarray} 
Therefore \begin{eqnarray} \label{e27a}
\!\!\!\! \|\Upsilon(\uz)-\Upsilon(\uz)\big\|_2  
&\le &   C_\l N^{1/2-\s} \td_2(N)  (\log\log \nu)^{1/2}  \Big( \sum_{j=1}^\nu    {|z_{j} - z'_{j}|^2\over p_j }\Big)^{1/2} 
    . 
\end{eqnarray} 
A Gaussian metric appears:   let indeed $g_1,\ldots, g_\nu$ be independent $ \mathcal N (0,1)$ distributed random
variables. Then
 $U(z):=   \sum_{j=1}^\nu  g_j p_j^{-1/2} z_{j} $ satisfies 
$$\|U(z)-U(z')\|_2=  \Big( \sum_{j=1}^\nu    {|z_{j} - z'_{j}|^2\over p_j }\Big)^{1/2} . 
$$

And so 
\begin{equation}\label{slep}\big\|\Upsilon(\uz)-\Upsilon(\uz)\big\|_2\le C_\l N^{1/2-\s} \td_2(N)  (\log\log \nu)^{1/2} \|U(z)-U(z')\|_2.
\end{equation}
Now  we  take again advantage of the comparison properties of Gaussian processes, and  deduce from Slepian's Lemma 
$$\E\,  \sup_{\uz,\uz'\in T^\nu} |\Upsilon (\uz')-\Upsilon(\uz)|\le C_\l N^{1/2-\s} \td_2(N) (\log\log \nu)^{1/2}\E\,  \sup_{\uz,\uz'\in
T^\nu} | U  (\uz')- U (\uz)|
 .  $$ 

But obviously 
$$ \sup_{\uz \in T^\nu}  |U(z)| =     \sum_{j=1}^\nu  |g_j |p_j^{-1/2}. $$
Thereby 
$$\E\,  \sup_{\uz'\in T^\nu} | U  (\uz')- U (\uz)|
\le C   \sum_{j=1}^\nu   p_j^{-1/2}\le C   \sum_{j=1}^\nu   {1\over (j\log j)^{1/2}} \le C \big({\nu \over \log \nu} \big)^{1/2}. $$
And by reporting 
$$\E\,  \sup_{\uz'\in T^\nu} |\Upsilon (\uz')-\Upsilon(\uz)|\le C_\l N^{1/2-\s} \td_2(N)  \big({\nu \log\log \nu\over \log \nu} \big)^{1/2}
 .  $$ 
  Observe also that
\begin{equation} \label{e210}
\|\Upsilon (\uz)\|_2\le C N^{1/2-\s} \td_2(N) ,\quad \uz\in
\T^{\nu}.
\end{equation}
Thus \begin{equation} \label{e211}
\E\,  \sup_{\uz'\in  T^\nu} |\Upsilon (\uz')|
\le    C
N^{1/2-\s} \td_2(N) \big({\nu \log\log \nu\over \log \nu} \big)^{1/2}.
\end{equation}
This   is   slightly better than in  \cite{LW2}, inequality (22), where one has the bound $ C
N^{1/2-\s} \td_2(N)  \nu^{1/2}$. By substituting in  (\ref{Q1A})      we  
get   
\medskip\par
\begin{equation} \label{Q_1}
   \E\, \sup_{\uz \in \T^{\nu }}\big|Q_1(\uz)\big|\le
  C_{\s,\l}  \,
N^{1/2-\s} \td_2(N) \big({\nu \log\log \nu\over \log \nu} \big)^{1/2},
\end{equation}
which is (\ref{Q1B}). 

\medskip\par

   Since $\td_1(N/p_\nu)\le \td_2(N/p_\nu)\le \td_2(N )$, we consequently get from (\ref{Q1B}),(\ref{Q2D}) and (\ref{e11}),   
\begin{equation}  \label{coll} \E\, \sup_{t \in \R} |\mathcal D (\s+it) |\le
  C_{\s,\l}  \,N^{1/2-\s}\td_2(N)  \, \bigg[ \big({\nu \log\log \nu\over \log \nu} \big)^{1/2}+
          { \tau^{1/2} \over (\log\tau)^{1/2}}+{   N^{1/2 }  \over \nu^{1/2} \log\nu} \bigg]
 .
\end{equation} 
  We now observe that $f(x):= (x\log\log x)^{1/2}(\log x)^{-1/2} + N^{1/2}x^{-1/2}(\log x)^{-1}$ satisfies 
 $$f'(x) \sim  \frac{1}{2}x^{-1/2}(\log x)^{-1/2}\big[ (\log\log x)^{1/2}-N^{1/2}x^{-1 }(\log x)^{-1/2} \big]. $$
Thus we choose 
$$\nu \sim {N^{1/2}\over (\log\log N)^{1/2}( \log N)^{1/2}}. $$
We get 
$${N^{1/2}\over \nu^{ 1/2} \log \nu  }\approx {N^{1/4}(\log\log N)^{1/4}\over (\log N)^{3/4}}\approx \Big({\nu\log\log \nu\over \log
\nu}\Big)^{1/2} .
$$

We find 
\begin{equation}  \label{coll1} \E\, \sup_{t \in \R} |\mathcal D (\s+it) |\le
   C_{\s,\l}  \,
N^{1/2-\s}\td_2(N)  \, \bigg[  {N^{1/4}(\log\log N)^{1/4}\over (\log N)^{3/4}}+
          { \tau^{1/2} \over (\log\tau)^{1/2}}  \bigg]
 .
\end{equation} 
  We also observe that ${N^{1/4}(\log\log N)^{1/4}\over (\log N)^{3/4}}\le  { \tau^{1/2} \over (\log\tau)^{1/2}}$, iff $\tau\ge (
{N\log\log N\over \log N} )^{1/2}$. Further when  $\tau\le ( {N\log\log
N\over \log N} )^{1/2}$,  we may also just set
$\nu=\tau$ in the initial decomposition, and thus ignore $Q^\e_2 $. It  means that we use the bound (\ref{Q_1}) in place of (\ref{coll}). This makes sense
when
$\tau$ is sufficiently small, namely when 
$ ({ \tau  \log \log \tau  \over  \log\tau  })^{1/2}\le {N^{1/4}(\log\log N)^{1/4}\over (\log N)^{3/4}}$; which is so when 
$\tau\le ( {N \over (\log N)\log\log
N} )^{1/2}$.  We
consequently have to distinguish three cases.
\medskip\par
{\bf Case 1.}\ $\big( {N\log\log N\over \log N}\big)^{1/2}  \le
\tau\le \pi(N).$
We get from  (\ref{coll1})
\begin{equation}  \label{coll11a} \E\, \sup_{t \in \R} |\mathcal D (\s+it) |\le
  C_{\s,\l}  \, {N^{1/2-\s}\td_2(N)  \tau^{1/2} \over (\log N)^{1/2}} .
\end{equation} 

{\bf Case 2.}\ $\big( {N \over (\log N)\log\log N}\big)^{1/2}\le
\tau \le \big( {N\log\log N\over \log N}\big)^{1/2}.$
In this case we  obtain from  (\ref{coll1})
\begin{equation}  \label{coll1b} \E\, \sup_{t \in \R} |\mathcal D (\s+it) |\le
  C_{\s,\l}  \,   {N^{3/4-\s}\td_2(N)(\log\log N)^{1/4}\over (\log N)^{3/4}} 
 .
\end{equation} 

{\bf Case 3.}\   $1\le \tau \le
\big( {N \over (\log N)\log\log N}\big)^{1/2}.$
By the comment made above, $\tau $ is small enough, and we forget   $Q^\e_2 $. We 
obtain from (\ref{Q_1}) directly
\begin{equation}  \label{coll1c} \E\, \sup_{t \in \R} |\mathcal D (\s+it) |\le
   C_{\s,\l}  \,
N^{1/2-\s}\td_2(N)  \,     
          \big({ \tau \log\log\tau \over  \log\tau  }\big)^{1/2}.   
 \end{equation} 
 
 \bigskip\par
Summarizing $$
 \E\, \sup_{t \in \R} |\mathcal D (\s+it) |\le  C_{\s,\l}  \, \td_2(N)\, B,$$
where $$
B=\cases{
  {N^{1/2-\s}   \tau^{1/2} \over (\log N)^{1/2}}
   &,\ if \ $\big( {N\log\log N\over \log N}\big)^{1/2}  \le
\tau\le \pi(N) ,$  \cr
& \cr
  {N^{3/4-\s} (\log\log N)^{1/4}\over (\log N)^{3/4}}
   &,\ if \  $\big( {N \over (\log N)\log\log N}\big)^{1/2}\le
\tau \le \big( {N\log\log N\over \log N}\big)^{1/2},$ \cr
& \cr
 N^{1/2-\s}        
          \big({ \tau \log\log\tau \over  \log\tau  }\big)^{1/2}
   &,\ if \  $1\le \tau \le
\big( {N \over (\log N)\log\log N}\big)^{1/2} .$
} 
 $$ This achieves the proof.  
 \cqfd
  
\medskip\par
 \noi
\section{Proof of Theorem \ref{t3}.}
  \medskip\par
 \noi  
   We   examine   more specifically the increments of the   Gaussian process   $\Upsilon$.  There is no loss   to assume 
$$p\mid K\quad \Rightarrow \quad p\le p_\nu. $$
  We have here
  \begin{eqnarray}
\Upsilon (\uz)    &=&
\sum_{n\in F_\nu    \atop (n,K)=1 }  n^{-\s}
\big\{\t_n \cos 2\pi \langle \ua(n),\uz\rangle +
\t_n'\sin 2\pi \langle \ua(n),\uz\rangle \big\}
 .
 \end{eqnarray}  
 And,  as  $(n,K)=1$ iff  $ a_\ell(n)>0\  \Rightarrow \  (p_\ell,K)=1 $,

  \begin{eqnarray*}
 \,\|\Upsilon(\uz)-\Upsilon(\uz)\big\|_2^2
&=& 4 \sum_{n\in F_\nu    \atop (n,K)=1 }    n^{-2\s}  \sin^2(\pi \langle\ua(n),\uz -\uz'\rangle)
 \cr
&\le & 4\pi^2  \sum_{n\in F_\nu    \atop (n,K)=1 }   n^{-2\s}
    \Big[ \sum_{     1\le j\le \nu\atop  (p_j,K)=1} a_j(n) |z_j - z'_j|\Big]^2
  .
\end{eqnarray*}
Now 
$$\sum_{n\in F_\nu    \atop (n,K)=1 }   n^{-2\s}
    \Big[ \sum_{     1\le j\le \nu\atop  (p_j,K)=1} a_j(n) |z_j - z'_j|\Big]^2=\sum_{n\in F_\nu    \atop (n,K)=1 } n^{-2\s}
    \sum_{     1\le j\le \nu\atop  (p_j,K)=1} a_j(n)^2 |z_j - z'_j| ^2$$$$ +\sum_{n\in F_\nu    \atop (n,K)=1 } n^{-2\s}
      \sum_{1\le j_1\not=j_2\le \nu\atop 
 (p_{j_1}p_{j_2},K)=1  } a_{j_1}(n)a_{j_2}(n) |z_{j_1} - z'_{j_1}|\ |z_{j_2} - z'_{j_2}|
 :=   S+R.$$
Further  
\begin{eqnarray}\label{debutR} R  &\le&  \  \sum_{1\le j_1\not=j_2\le \nu\atop 
 (p_{j_1}p_{j_2},K)=1  }  |z_{j_1} - z'_{j_1}|\ |z_{j_2} -z'_{j_2}|
    \sum_{b_1,b_2=1}^\infty b_1 b_2 \sum_{n\in F_\nu, \, (n,K)=1 \atop  a_{j_1}(n)=b_1,\, a_{j_2}(n)=b_2 } {1\over
n^{ 2\s}}
 \cr&\le &  C \sum_{1\le j_1\not=j_2\le \nu\atop 
 (p_{j_1}p_{j_2},K)=1  } |z_{j_1} - z'_{j_1}|  |z_{j_2} - z'_{j_2}|
    \sum_{b_1,b_2=1}^\infty      { b_1 b_2 \over  (p_{j_1}^{  b_1 }     p_{j_2}^{  b_2 })^{2\s}
}
 \cr& &  \qquad\times
\Big[\sum_{m\le N/(p_{j_1}^{  b_1 }     p_{j_2}^{  b_2 }) }   m^{
-2\s}   \Big]
\cr&\le &  C N^{1-2\s}\sum_{1\le j_1\not=j_2\le \nu\atop 
 (p_{j_1}p_{j_2},K)=1  } |z_{j_1} - z'_{j_1}|  |z_{j_2} - z'_{j_2}|
    \sum_{b_1,b_2=1}^\infty      { b_1 b_2 \over   p_{j_1}^{  b_1 }     p_{j_2}^{  b_2 }  
}
  .
\end{eqnarray}
But
 $$\sum_{b=1}^\infty  {b   \over
p_k^{b } } =\sum_{b=1}^\infty  {b   \over
2^{b } } \Big[{2   \over
p_k  }\Big]^{ b }\le  {2   \over
p_k  } \sum_{b=1}^\infty  {b   \over
2^{b } } \le Cp_k^{-1  } . $$

Thus 
 \begin{eqnarray}\label{debutR} R   &\le &  C N^{1-2\s}\Big(\sum_{1\le j \le \nu\atop 
 (p_{j } ,K)=1  } {|z_{j } - z'_{j }|\over p_j}  
    \Big)^2\cr &\le& C N^{1-2\s}\Big(\sum_{1\le j \le \nu\atop 
 (p_{j } ,K)=1  } {1\over p_j}  
    \Big) \Big(\sum_{1\le j \le \nu\atop 
 (p_{j } ,K)=1  } {|z_{j } - z'_{j }|^2\over p_j}  
    \Big) 
  .
\end{eqnarray} And
  \begin{eqnarray}\label{debutS}S&\le& \sum_{n\in F_\nu    \atop (n,K)=1 } n^{-2\s}
    \sum_{     1\le j\le \nu\atop  (p_j,K)=1} a_j(n)^2 |z_j - z'_j| ^2 \cr& \le &\sum_{     1\le j\le \nu\atop  (p_j,K)=1}|z_j - z'_j|
^2\sum_{b=1}^\infty b^2\sum_{{n\in F_\nu   
\atop (n,K)=1 }\atop a_j(n)=b} {1\over n^{2\s}} 
\cr& \le &\sum_{     1\le j\le \nu\atop  (p_j,K)=1}|z_j - z'_j|
^2\sum_{b=1}^\infty {b^2\over p_j^{2b\s}}\sum_{ m\le N/p_j^b   
   } {1\over m^{2\s}} 
\cr& \le &\sum_{     1\le j\le \nu\atop  (p_j,K)=1}|z_j - z'_j|
^2\sum_{b=1}^\infty {b^2\over p_j^{ b }}\le \sum_{     1\le j\le \nu\atop  (p_j,K)=1}   {|z_j - z'_j|
^2\over p_j }  .
 \end{eqnarray}
     
    Therefore,  
\begin{equation}  \label{estK1} \|\Upsilon(\uz)-\Upsilon(\uz)\big\|_2^2
 \le    C_\s  N^{1-2\s} \,      
  \Big[ \sum_{k\le \nu\atop
 p_k\not|K  }
  {|z_k - z'_k|^2\over   p_k }\Big]\max\Big(1,\sum_{1\le j \le \nu\atop 
 (p_{j } ,K)=1  } {1\over p_j}  
    \Big)
.
\end{equation}
Let 
$$ \D:=  N^{1/2- \s}  \,      
   \max\Big(1,\sum_{1\le j \le \nu\atop 
 (p_{j } ,K)=1  } {1\over p_j}  
    \Big)^{1/2} .$$
 We   obtain 
\begin{equation}  \label{estK2} \|\Upsilon(\uz)-\Upsilon(\uz)\big\|_2 
  \le   C_\s \D    
   \Big[ \sum_{k\le \nu\atop
 p_k\not|K  }
  {|z_k - z'_k|^2\over   p_k }\Big]^{1/2} 
.
\end{equation}  Let $g_1,\ldots, g_\nu$ be independent $ \mathcal N (0,1)$ distributed random
variables and define
 $U(z):= \sum_{k\le \nu\atop
 p_k\not|K  }  g_k p_k^{-1/2} z_{k} $. Then
\begin{equation}\label{slep}\big\|\Upsilon(\uz)-\Upsilon(\uz)\big\|_2\le C_\s \D \|U(z)-U(z')\|_2.
\end{equation}
  We deduce from Slepian's Lemma   
$$\E\,  \sup_{\uz'\in T^\nu} |\Upsilon (\uz')-\Upsilon(\uz)|\le C_\s \D \E\,  \sup_{\uz'\in T^\nu} | U  (\uz')- U (\uz)|
 .  $$ 
Obviously 
$$ \sup_{\uz \in T^\nu}  |U(z)| =   \sum_{k\le \nu\atop
 p_k\not|K  }  {|g_k |\over p_k^{ 1/2}}. $$
Thereby 
$$\E\,  \sup_{\uz'\in T^\nu} | U  (\uz')- U (\uz)|
\le C    \sum_{k\le \nu\atop
 p_k\not|K  }   p_k^{ -1/2} . $$
And by reporting 
$$\E\,  \sup_{\uz'\in T^\nu} |\Upsilon (\uz')-\Upsilon(\uz)|\le  C_\s \D\Big[\sum_{k\le \nu\atop
 p_k\not|K  } {1\over \sqrt {p_k} }\Big] 
 .  $$ 
  But
\begin{equation} \label{e210}
\|\Upsilon (\uz)\|_2 \le \Big[ \sum_{n\in F_\nu\atop (n,K)=1}{1\over n^{ 2\s}}\Big]^{1/2}\le C_\s N^{1/2-\s}\Big[ 
\sum_{k\le \nu\atop
 p_k\not|K  }{1\over
p_j}\Big]^{1/2} ,\qquad \uz\in
\T^{\nu}.
\end{equation}
Thus \begin{equation} \label{e212}
\E\,  \sup_{\uz'\in  T^\nu} |\Upsilon (\uz')|
\le   C_\s \D\Big[\sum_{k\le \nu\atop
 p_k\not|K  } {1\over \sqrt {p_k} }\Big] ,
\end{equation}
 or \begin{equation} \label{e212}
\E\,  \sup_{\uz'\in  T^\nu} |\Upsilon (\uz')|
\le  N^{1/2- \s}  \,      
   \max\Big(1,\sum_{k\le \nu\atop
 p_k\not|K  } {1\over   {p_k} }   
    \Big)^{1/2}\Big[\sum_{k\le \nu\atop
 p_k\not|K  } {1\over \sqrt {p_k} }\Big] .
\end{equation}
\cqfd
 

\medskip\par
 \noi
\section{Intermediate results.}
  \medskip\par
 \noi

    The following   result of Hall   will be useful.
 Let $f$ be defined on positive integers and satisfying $f(1)=1$, $0\le f(n)\le 1$, and being sub-multiplicative.
 \smallskip\par
Put 
$$  {\bf\Pi}_x(f)= \prod_{p\le x}\Big(1-{1\over p}\Big)\Big(1+ {f(p)\over p}+ {f (p^2)\over p^2}+ \ldots\Big)
$$
Then (\cite{Hal}, theorem 2)
 \begin{equation}\label{hal3}
 \sum_{n\le x} f(n)\le C\ x {\bf\Pi}_x(f),
\end{equation}
$C$ being an absolute constant. This estimate allows in turn a similar control   for bounded non-negative sub-multiplicative functions. 
\smallskip \par  Apply it to $f=d_K$. As  $1+ {f(p)\over p}+ {f (p^2)\over p^2}+
\ldots = {p \over p-1}$, if
$(p,K)=1$, we have 
\begin{equation}  {\bf\Pi}_x(f) = \prod_{p\le x\atop (p,K)>1}\big(1-{1\over p}\big)= 
 \prod_{p\le x\atop p|K}\big(1-{1\over p}\big) .
\end{equation}
 
Hence the  classical estimate, (see  \cite{Hal} for references) 
  \begin{equation}\label{hal3}
 \p_K(x):=\#\big\{ k\le x: (k,K)=1\big\} \le C\ x
 \prod_{  p|K \atop p\le x  }\big(1-{1\over p}\big).
\end{equation}

    We will need the following technical Lemma.
\begin{lem}\label{lem3} a) Let a real $\b>0$ and integer $L>0$. Then   
\begin{equation}\label{lem3e1}\sum_{(n,L)=1\atop n\le x} n^{-\b}  
\le C_\b x^{1-\b}\prod_{p |L\atop p\le x }\big(1-{1\over p}\big)  . 
\end{equation} 
b) Let  $0\le \b<1$. Then   
  \begin{equation}\label{pplus} \sum_{ n\le
x \atop P^+(n)\le y  }  {1   \over n^{  
\b}}\le  C_\b x^{   1-\b}   e^{-{1\over 2}{\log x\over \log y}} ,
\end{equation} 
for some constant  $C_\b $, 
    $y\ge y_\b$,     $y/x\le c_\b  $.

\smallskip\par \noi c)  If $\b=1$, then 
   \begin{eqnarray} \sum_{y\le  n\le
x \atop P^+(n)\le y  }  {1   \over n } \le C\log y .
 \end{eqnarray} \end{lem}

\begin{rem}\rm It is natural to compare, in our setting, estimates a) anf b), via the relation
$$\sum_{   P^+(n)\le p_\tau \atop n\le
N   }  {1   \over n^{  
\b}}= \sum_{(n,K_\tau)=1\atop n\le N} n^{-\b}$$
where $K_\tau$ is defined in (\ref{ktau}). By a) and Mertens Theorem we get 
$$ \sum_{   P^+(n)\le p_\tau \atop n\le
N   }  {1   \over n^{  
\b}} 
\le C_\b N^{1-\b}\prod_{\tau<\ell\le \pi(N) }\big(1-{1\over p}\big)\le C_\b N^{1-\b}{\log p_\tau\over \log N}  . 
$$
      However, by using b)  we get the much better bound $  C_\b N^{   1-\b}   e^{-{1\over 2}{\log N\over \log p_\nu}}$. 
\end{rem} 
\medskip\par \noi 
{\it Proof.}  a) By applying formula (\ref{simple}) with $a_m=\chi\{(m,L)=1\}$,     $b_m=m^{-\s}$, $1\le m\le x$,
$$\sum_{(m,L)=1\atop m\le x} m^{-\b}\le {A(x)\over x^\b} +\b\int_1^x A(t)  {dt\over t^{\b+1}} ,$$ 
where $A(t)=\sum_{n<t} d_L(n)$. 
\smallskip\par
But by   Hall's estimate (\ref{hal3}), $A(t)\le Ct \prod_{p|L\atop p\le t}\big(1-{1\over p}\big)$. 
 Thus
$$\sum_{(m,L)=1\atop m\le x} m^{-\b}\le Cx\prod_{p|L\atop p\le x}\big(1-{1\over p}\big)\,{1\over x^\b} +C\b\int_1^x \prod_{p|L\atop p\le
t}\big(1-{1\over p}\big)  {dt\over t^{\b }} \le C_\b\int_1^x \prod_{p|L\atop p\le
t}\big(1-{1\over p}\big)  {dt\over t^{\b }} .$$ 
Applying now twice Mertens's theorem, gives
\begin{eqnarray}\prod_{p|L\atop p\le t}\big(1-{1\over p}\big)&=&{\prod_{  p\le t}\big(1-{1\over p}\big)\over \prod_{p\not|L\atop p\le t}\big(1-{1\over
p}\big)}
\le {C\over \log t\, \prod_{p\not|L\atop p\le x }\big(1-{1\over p}\big)} \le
{C\prod_{p |L\atop p\le x }\big(1-{1\over p}\big)\over \log t\, \prod_{p \le x }\big(1-{1\over p}\big)} \cr & \le &C\,
{\log  x \over \log  t}\prod_{p |L\atop p\le x }\big(1-{1\over p}\big). 
\end{eqnarray}
Hence
$$\sum_{(m,L)=1\atop m\le x} m^{-\b}  \le C_\b\log x \prod_{p |L\atop p\le x }\big(1-{1\over p}\big)\int_1^x    {dt\over t^{\b }\log  t} 
\le C_\b x^{1-\b}\prod_{p |L\atop p\le x }\big(1-{1\over p}\big)  .$$ 
\smallskip\par b)   Let $\Psi(x,y):=\#\{n\le x: P^+(n)\le y\}$.  
 By using this time (\ref{simple}) with   $a_n=\chi\{P^+(n)\le y\}$   $1\le n\le N$,  we obtain
\begin{eqnarray} \sum_{1\le  n\le
x \atop P^+(n)\le y  }  {1   \over n^{  
\b}}&= &{\#\{1\le  n\le
x:P^+(n)\le y\}\over x^\b}  +\b\int_1^x{\#\{1\le  n\le
t:P^+(n)\le y\}\over t^{\b+1}}dt\cr &= &{\Psi(x,y)\over x^\b}+\b\int_1^y{dt\over t^{\b }}\, +\b\int_y^x{\Psi(t,y)\over t^{\b+1}}\,dt.
\end{eqnarray}
  Recall that $\Psi(x,y)\le xe^{-{1\over 2}{\log x\over \log y}}  $, $x\ge y\ge 2$,
(\cite{T}, Chapter III.5).  Thus, for $y$ sufficiently large    to have $1   -\b>{1\over  \log
y}$,  
\begin{eqnarray}\int_y^x{\Psi(t,y)\over t^{\b+1}}\,dt&\le& \int_y^xe^{-{1\over 2}{\log t\over \log y}}{dt\over t^{\b }}=\int_y^xt^{- {1\over 2\log y}-\b}\,
dt  
\cr &=&{ 1\over 1-{1\over 2\log y}-\b}\ \Big(t^{1- {1\over 2\log y}-\b}\,\Big|_{t=y}^{t=x}\le { 2\over 1- \b}\  x^{1- {1\over 2\log
y}-\b} \cr &=& { 2\over 1- \b}\ x^{   1-\b}   e^{-{1\over 2}{\log x\over \log y}}  .\end{eqnarray}
Therefore \begin{eqnarray} \sum_{y\le  n\le
x \atop P^+(n)\le y  }  {1   \over n^{  
\b}}&\le & C_\b\Big[ x^{   1-\b}   e^{-{1\over 2}{\log x\over \log y}}+  y^{1-\b} \Big] .
 \end{eqnarray}
 Now, we have $x^{   1-\b}   e^{-{1\over 2}{\log x\over \log y}}\ge   y^{1-\b} $ iff $\log {x\over y} \ge {1\over 2(1-\b)}\, {\log x\over \log y}$. Write
$x=\theta y$, $\theta\ge 1$. This means
$$\log \theta \ge {1\over 2(1-\b)}\, {\log \theta y\over \log y}={1\over 2(1-\b)}\, \big\{ {\log \theta  \over \log y} +1 \big\}, $$
or 
$$\log \theta \, \Big\{1- {1  \over 2(1-\b) \log y}   \Big\}\ge  {1\over 2(1-\b)}. $$
 If $y$ is large enough, $y\ge y_\b$,   $y/x$ small enough, $y\le c_\b x$, then the above condition is satisfied. Consequently
  \begin{equation}\label{pplus} \sum_{ n\le
x \atop P^+(n)\le y  }  {1   \over n^{  
\b}}\le  C_\b x^{   1-\b}   e^{-{1\over 2}{\log x\over \log y}}  .
\end{equation}

 c) The case $\b=1$ can be treated as before:
 \begin{eqnarray} \sum_{1\le  n\le
x \atop P^+(n)\le y  }  {1   \over n }  &= &{\Psi(x,y)\over x }+ \int_1^y{dt\over t }\, + \int_y^x{\Psi(t,y)\over t }\,dt.
\end{eqnarray}
 And 
\begin{eqnarray}\int_y^x{\Psi(t,y)\over t^2}\,dt&\le& \int_y^xe^{-{1\over 2}{\log t\over \log y}}{dt\over t }=\int_y^xt^{- {1\over
2\log y}-1}\, dt  
={ 1\over  -{1\over 2\log y} }\ \Big[t^{ - {1\over 2\log y} }\,\Big|_{t=y}^{t=x}\cr&\le& { 1\over   {1\over 2\log y} }\  y^{ - {1\over
2\log y} }  \le   C  \log y   .\end{eqnarray}
Therefore \begin{eqnarray} \sum_{y\le  n\le
x \atop P^+(n)\le y  }  {1   \over n }&\le & C \Big[     e^{-{1\over 2}{\log x\over \log y}}+   \log y  \Big]\le C\log y .
 \end{eqnarray}
   One can  however  get this directly.   Let $ j=j_y= \max\{\ell:p_\ell\le y\}$. Then,  for any $\b>0$,
 \begin{equation}\label{direct} \sum_{1\le  n\le x \atop P^+(n)\le y  }  {1   \over n^\b }\le
\sum_{\a_1=0}^\infty\ldots\sum_{\a_j=0}^\infty {1\over p_1^{\a_1\b}\ldots p_j^{\a_j\b}}=\prod_{\ell=1}^j \Big({1\over 1-{1\over
p_\ell^\b}}\Big) .
\end{equation}
 And   when $\b=1$, by    Mertens Theorem, the latter is less than $\le C\log y$.\cqfd   
 \smallskip \par

This last argument can serve to get a two-sided estimate when $y$ is not too large. In this case, the estimates depend on $y$ only.
\begin{lem} \label{util}If $y=o(\log x)$, then we have
  for any $\b>0$,
 \begin{equation}\label{direct3}c_\b \prod_{\ell =1}^j\Big[{1\over 1-{1\over p_\ell^\b}} \Big]\le  \sum_{1\le  n\le x \atop P^+(n)\le y 
}  {1\over n^\b}\le C_\b
\prod_{\ell =1}^j\Big[{1\over 1-{1\over p_\ell^\b}} \Big].
\end{equation}
And the involved constants $c_\b, C_\b$ depend on $\b$ only. In particular 
\begin{equation}\label{direct2}  C_1\log y\le \sum_{1\le  n\le x \atop P^+(n)\le y  }  {1\over n} \le C_2\log y.
\end{equation}
 \end{lem}
 \smallskip\par\noi  {\it Proof.} Indeed,  notice first, as   $p_j\sim j\log j$, that we have $j\le Cy/\log y$. Now consider integers   $n=
p_1^{\a_1}\ldots p_j^{\a_j}$, such  that 
 $\max
\{\a_\ell,
\ell\le j\}\le  H:= ( {\log x } )/Cy$.    Thus
 $$n\le y^{j\max \{\a_\ell, \ell\le j\}} =e^{j(\log y)\max \{\a_\ell, \ell\le j\}}\le e^{ {Cy\over \log y}(\log y)  \{{\log x\over
Cy}\}}\le x.$$
 We may also assume that $(H+1) \b\ge 2$.  Therefore
 $$\sum_{1\le  n\le x \atop P^+(n)\le y  }  {1\over n^\b}\ge \sum_{\a_1=0}^H\ldots\sum_{\a_j=0}^H
{1\over p_1^{\a_1\b}\ldots p_j^{\a_j\b}} = \prod_{\ell =1}^j\Big[{1\over 1-{1\over p_\ell^\b}}- \sum_{\a_j=H+1}^\infty 
{1\over p_\ell^{  \a_\ell\b}}\Big]$$
$$= \prod_{\ell =1}^j\Big[{1\over 1-{1\over p_\ell^\b}} \Big]\prod_{\ell =1}^j\Big[1- 
{1\over p_\ell^{(H+1)\b}}\Big]\ge c_\b \prod_{\ell =1}^j\Big[{1\over 1-{1\over p_\ell^\b}} \Big] .$$
But
$$\prod_{\ell =1}^j\Big[1- 
{1\over p_\ell^{(H+1)\b}}\Big] \ge\prod_{\ell =1}^j\Big[1- 
{1\over p_\ell^{2 }}\Big] \ge      e^{-C' \sum_{\ell =1}^\infty
p_\ell^{-2}}   >0.
$$
since the series $\sum_{\ell =1}^\infty
p_\ell^{-2} $ is obviously convergent.  And so, in view of (\ref{direct}) 
 \begin{equation}\label{direct2}c_\b \prod_{\ell =1}^j\Big[{1\over 1-{1\over p_\ell^\b}} \Big]\le  \sum_{1\le  n\le x \atop P^+(n)\le y 
}  {1\over n^\b}\le C_\b
\prod_{\ell =1}^j\Big[{1\over 1-{1\over p_\ell^\b}} \Big].
\end{equation}
When $\b=1$,  by using Mertens Theorem
$$C_1\log y\le \sum_{1\le  n\le x \atop P^+(n)\le y  }  {1\over n}\le C_2\log y.$$\cqfd
 We continue with some other  useful observations.   
\begin{rem}\rm Let   $u:={\log x\over \log y}$ and   $\rho(.)$ denote Dickman's function. According to (\cite{T}, p.435),  
\begin{eqnarray}\sum_{   n\le x
  \atop P^+(n)\le y  }  {1   \over n   
 } &=&\log y\int_0^u\rho(v) dv + \mathcal O(u)=\log y\Big( e^\gamma+ \mathcal O\big( {u\over \log y} + e^{-u/2}\big)\Big)+
\mathcal O(u)\cr &= & e^\gamma \log y +\mathcal O(u),
\end{eqnarray}
for $x\ge y\ge 2$,   $\gamma$ being Euler's constant. 
 \end{rem}   
\bigskip\par   In \cite{LW1}, we introduced a new approach to lower bounds. It will be necessary to briefly recall  its principle.  We  begin with the   
lemma below (\cite{LW1}, Lemma 3.1).
 \begin{lem} Let  $X=\{X_z, z\in Z\}$ and $Y=\{Y_z,
z\in Z\}$ be two finite sets of random variables defined on a
common probability space. We assume that $X$ and $Y$ are
independent and that the random variables $Y_z$ are all centered.
Then 
$$ \E\sup_{z\in Z}|X_z + Y_z| \ge  \E\sup_{z\in Z}|X_z  |.$$
\end{lem} 
 
 Let $\ud =\{d_n, n\ge 1\}$ be a   sequence of reals. By the reduction step (\ref{e11})
$$\sup_{t\in \R}\big|\sum_{j=1}^\tau \sum_{n\in E_j} d_n\e_n
n^{ -\s - it}\big| =\sup_{\uz \in \T^\tau}\big|Q(\uz)\big|.
$$
  where
 $$ Q(\uz)= \sum_{j=1}^\tau \sum_{n\in E_j}
   d_n \e_n n^{-\s} e^{2i\pi\langle \ua(n),\uz\rangle}.
$$
  Introduce the following subset  of $\T^\tau$,  
$$\zz=\Big\{ \uz=\{z_j,  1\le j\le \tau\}    :   \hbox{$z_j=0$,   
if $j\le \tau/2$,\  and
\  $z_j\in\{0,1/2\}$, if $j\in]\tau/2,\tau]$} \Big\} .
$$
 Observe that for any $\uz\in \zz$, any $n$, $ e^{2i\pi\langle \ua(n),\uz\rangle} = \cos(2\pi\langle
\ua(n),\uz\rangle) = (-1)^{2\langle \ua(n),\uz\rangle}$. It follows that   $\Im\, Q(\uz)=0$, and so  
$$
 Q(\uz)
= \sum_{\tau/2<j\le \tau}  \sum_{n\in E_j}
d_n\e_n n^{-\s} {(-1)}^{2\langle \ua(n),\uz\rangle },
 \qq \uz\in \zz. $$
 Thereby the restriction of $Q$   to $\zz$ is just a finite
 rank Rademacher process.
 Now define
 $$ {\cal L}_j=\Big\{n=p_j \, \tilde n\ : \
  \tilde n\le {N\over p_{j}}\ \hbox{and}\
P^+(\tilde n)\le p_{\tau/2}\Big\}, \qq  \qq j\in(\tau/2,\tau]. $$
 Since
  $  E_j \supset{\cal L}_j$,
 $j=1,\ldots \tau,$
 the sets ${\cal L}_j$ are pairwise disjoint.
 Put for $z\in \zz$,
$$
 Q'(\uz)
= \sum_{  \tau/2<j\le \tau }\sum_{n\in  {\cal L}_j}
\e_n n^{-\s} {(-1)}^{2\langle \ua(n),\uz\rangle }
 .
 $$
 Since $\{Q(\uz)-Q'(\uz),\uz\in \zz\}$ and $\{
Q'(\uz),\uz\in \zz\}$ are independent, we deduce from the above Lemma that
$$\E \sup_{\uz\in \zz} |Q(\uz)|
 \ge
  \E \sup_{\uz\in \zz} \left|Q'(\uz) \right|.
$$

It is possible to  proceed to a direct evaluation of $Q'(\uz)$ and we recall that
$$\sup_{\uz\in \zz} \left|Q'(\uz) \right|
= \sum_{  \tau/2<j\le \tau }\big|\sum_{n\in  {\cal L}_j}
d_n\e_n  n^{-\s} \big|, 
$$ 
which, in view of the Khintchine inequalities for Rademacher sums, allows to get (\cite{LW1}, Proposition 3.2)\begin{prop}
 There exists a universal
constant $c$ such that for any system of coefficients $(d_n)$ 
$$c\ \sum_{  \tau/2<j\le \tau }
\big|\sum_{n\in  {\cal L}_j}
d_n^2  n^{-2\s} \big|^{1/2} \le \E\, \sup_{\uz\in \zz}
\left|Q'(\uz) \right|
\le
\sum_{  \tau/2<j\le \tau }
\big|\sum_{n\in  {\cal L}_j} d_n^2  n^{-2\s} \big|^{1/2}.
$$
\end{prop}\medskip

 Consequently
\begin{equation}\label{minor1}\E \sup_{t\in \R}\big|\sum_{j=1}^\tau \sum_{n\in E_j} d_n\e_nn^{-\s-it}\big|\ge c\ \sum_{  \tau/2<j\le \tau }
\big|\sum_{n\in  {\cal L}_j}
d_n^2  n^{-2\s} \big|^{1/2}. 
\end{equation} 

 \medskip\par 
\section{Proof of Theorem \ref{smalltau}.}
\medskip\par\noi 

\noi {\it Proof of the lower bound.} Take $d_n\equiv
1$   in estimate (\ref{minor1}).  We get
 \begin{eqnarray}\E \sup_{t\in \R}\big|\sum_{  n\le N\atop P^+(n)\le p_{\tau}}  {\e_n
\over n^{ \s+ it}}\big|
 &=&\E \sup_{t\in \R}\big|\sum_{j=1}^\tau \sum_{n\in E_j}  {\e_n\over n^{\s+it}} \big|\cr  & \ge&c\ \sum_{  \tau/2<j\le \tau }
\big|\sum_{n\in  {\cal L}_j}
   {1\over n^{2\s}} \big|^{1/2}.
\end{eqnarray}
By   assumption
 $   \log{N\over p_{j}}\ge \log {N\over p_{\tau/2} }=\log {N }-\log {  p_{\tau/2} }\gg  p_{\tau/2}  $.  
 Owing to the very definition of the sets ${\cal L}_j$, and using  Lemma \ref{util}, we get 
 \begin{eqnarray}\sum_{  \tau/2<j\le \tau }
\big|\sum_{n\in  {\cal L}_j}
   {1\over n^{2\s}} \big|& =&\sum_{  \tau/2<j\le \tau }
{1\over   p_j^\s}\Big[\sum_{\tilde n\le {N\over p_{j}}\atop P^+(\tilde n)\le p_{\tau/2}}
   {1\over \tilde n^{2\s}} \Big]^{1/2}\cr& \ge & C_\s \prod_{\ell =1}^\tau\Big[{1\over 1-{1\over p_\ell^{2\s}}} \Big]^{1/2}\sum_{  \tau/2<j\le \tau }
{1\over p_j^\s} 
\cr& \ge & C_\s \prod_{\ell =1}^\tau \Big[{1\over 1-{1\over p_\ell^{2\s}}} \Big]^{1/2}   {\tau^{1-\s}  \over( \log \tau)^\s }  .
\end{eqnarray}
   Consequently
\begin{equation}\label{lowersmalltau}
\E \sup_{t\in \R}\big|\sum_{  n\le N\atop P^+(n)\le p_{\tau}}  {\e_n
\over n^{ \s+ it}}\big|
 \ge C_\s \prod_{\ell =1}^\tau \Big[{1\over 1-{1\over p_\ell^{2\s}}} \Big]^{1/2}  {\tau^{1-\s}  \over( \log \tau)^\s } .
\end{equation}
And if $\s=1/2$, by Mertens Theorem,
 \begin{equation}\label{lowerdemismalltau}\E \sup_{t\in \R}\big|\sum_{  n\le N\atop P^+(n)\le p_{\tau}}  {\e_n
\over n^{  {1\over 2}+ it}}\big|
\ge    C    {\tau }^{1/2}.\end{equation}  
     
\medskip\par
 \noi
 {\it Proof of the upper bound.}
    We have   \begin{eqnarray*}
\Upsilon (\uz)    &=&
\sum_{n\in F_\tau  } {1\over   n^\s}
\big\{\t_n \cos 2\pi \langle \ua(n),\uz\rangle +
\t_n'\sin 2\pi \langle \ua(n),\uz\rangle \big\}
 .
 \end{eqnarray*}  
 And   
  $\   \|\Upsilon(\uz)-\Upsilon(\uz')\big\|_2^2
  \le   4\pi^2  \sum_{n\in F_\tau    }   {1\over  n^{2\s}}
    \big[ \sum_{j=1}^\tau a_j(n) |z_j - z'_j|\big]^2
  $.
 \smallskip\par
Now 
$$\sum_{n\in F_\tau     }  {1\over  n^{2\s}}
    \Big[ \sum_{j=1}^\tau a_j(n) |z_j - z'_j|\Big]^2=\sum_{n\in F_\tau     }
  {1\over  n^{2\s}}    \sum_{j=1}^\tau a_j(n)^2 |z_j - z'_j| ^2$$$$ +\sum_{n\in F_\tau    }   {1\over  n^{2\s}}
      \sum_{1\le j_1\not=j_2\le \tau   } a_{j_1}(n)a_{j_2}(n) |z_{j_1} - z'_{j_1}|\ |z_{j_2} - z'_{j_2}|
 :=   S+R.$$
Further, by using Lemma \ref{util}  
   
\begin{eqnarray*}\label{debutR} R  &\le&  \  \sum_{1\le j_1\not=j_2\le \tau   }  |z_{j_1} - z'_{j_1}|\ |z_{j_2} -z'_{j_2}|
    \sum_{b_1,b_2=1}^\infty b_1 b_2 \sum_{{n\in F_\tau    \atop  a_{j_1}(n)=b_1} \atop a_{j_2}(n)=b_2 }  {1\over  n^{2\s}} \cr&\le &  C \sum_{1\le j_1\not=j_2\le \tau   } |z_{j_1} - z'_{j_1}|  |z_{j_2} - z'_{j_2}|
    \sum_{b_1,b_2=1}^\infty      { b_1 b_2 \over   p_{j_1}^{ 2 b_1\s }     p_{j_2}^{ 2 b_2 \s}  
}
 \Big[\sum_{m\le N/(p_{j_1}^{  b_1 }     p_{j_2}^{  b_2 })\atop P^+(m)\le p_\tau  }   {1\over  m^{2\s}} \Big]
\cr&\le &  C_\s \prod_{\ell =1}^\tau\Big[{1\over 1-{1\over p_\ell^{2\s}}} \Big]\sum_{1\le j_1\not=j_2\le \tau   } {|z_{j_1} - z'_{j_1}|\over p_{j_1}^{2\s}}{  |z_{j_2} - z'_{j_2}|\over p_{j_2}^{2\s}}
   .
\end{eqnarray*}
 Thus 
 \begin{eqnarray}\label{debutR1} R   & \le &   C_\s \prod_{\ell =1}^\tau\Big[{1\over 1-{1\over p_\ell^{2\s}}} \Big]\Big( \sum_{j=1}^\tau {|z_{j } - z'_{j }|\over p_j^{2\s}}  
    \Big)^2\cr & \le &  C_\s \prod_{\ell =1}^\tau \Big[{1\over 1-{1\over p_\ell^{2\s}}} \Big] \Big(\sum_{j=1}^\tau {1 \over p_j^{2\s}}  
    \Big) 
 \Big(\sum_{j=1}^\tau {|z_{j } - z'_{j }|^2\over p_j^{2\s}}  
    \Big) 
\cr & \le &  C_\s \prod_{\ell =1}^\tau \Big[{1\over 1-{1\over p_\ell^{2\s}}} \Big] \Big({ \tau^{1-2\s} \over(\log \tau)^{2\s} }
    \Big) 
 \Big(\sum_{j=1}^\tau {|z_{j } - z'_{j }|^2\over p_j^{2\s}}  
    \Big) 
  .
\end{eqnarray} And
  \begin{eqnarray}\label{debutS1}S&\le& \sum_{n\in F_\tau     }{1\over n^{2\s}}
    \sum_{j=1}^\tau a_j(n)^2 |z_j - z'_j| ^2   \le  \sum_{j=1}^\tau|z_j - z'_j|
^2\sum_{b=1}^\infty b^2\sum_{{n\in F_\tau    
  }\atop a_j(n)=b} {1\over n^{2\s} } 
\cr& \le &\sum_{j=1}^\tau|z_j - z'_j|
^2\sum_{b=1}^\infty {b^2\over p_j^{ 2b\s }}\sum_{ m\le N/p_j^b   
  \atop P^+(m)\le p_\tau  } {1\over m^{2\s} } 
\cr& \le & C_\s \prod_{\ell =1}^\tau \Big[{1\over 1-{1\over p_\ell^{2\s}}} \Big]\sum_{j=1}^\tau|z_j - z'_j|
^2\sum_{b=1}^\infty {b^2\over p_j^{ 2b\s }}
\cr& \le & C_\s \prod_{\ell =1}^\tau \Big[{1\over 1-{1\over p_\ell^{2\s}}} \Big]\sum_{j=1}^\tau   {|z_j - z'_j|
^2\over p_j^{2\s} }  .
 \end{eqnarray}
     
     Consequently
\begin{equation}  \label{estK1} \|\Upsilon(\uz)-\Upsilon(\uz)\big\|_2^2
 \le   C_\s \prod_{\ell =1}^\tau \Big[{1\over 1-{1\over p_\ell^{2\s}}} \Big]  \Big({ \tau^{1-2\s} \over(\log \tau)^{2\s} }
    \Big) 
 \,      
  \Big[ \sum_{j=1}^\tau
  {|z_j - z'_j|^2\over   p_j^{2\s} }\Big]\ 
.
\end{equation}
     We deduce from Slepian's Lemma,   noting that $\log p_\tau\sim \log \tau$   
 \begin{eqnarray*}\E\,  \sup_{\uz, \uz'\in T^\tau } |\Upsilon (\uz')-\Upsilon(\uz)|&\le &  C_\s \prod_{\ell =1}^\tau \Big[{1\over 1-{1\over p_\ell^{2\s}}} \Big]^{1/2} \Big({ \tau^{{1\over 2}- \s} \over(\log \tau)^{ \s} }
    \Big) 
 \Big[\sum_{j=1}^\tau {1\over
p_j^{ \s} }\Big] 
\cr&\le &  C_\s \prod_{\ell =1}^\tau \Big[{1\over 1-{1\over p_\ell^{2\s}}} \Big]^{1/2} \Big({ \tau^{{1\over 2}- \s} \over(\log \tau)^{ \s} }
    \Big) 
 \Big({ \tau^{1-  \s} \over(\log \tau)^{  \s} }
    \Big) 
\cr&= &  C_\s \prod_{\ell =1}^\tau \Big[{1\over 1-{1\over p_\ell^{2\s}}} \Big]^{1/2} \Big({ \tau^{{3\over 2}- 2\s} \over(\log \tau)^{2 \s} }
    \Big) .
    \end{eqnarray*} 
 But
\begin{equation} \label{e210}
\|\Upsilon (\uz)\|_2 \le \Big[ \sum_{n\in F_\tau  }{1\over n^{2\s} }\Big]^{1/2}\le C_\s \prod_{\ell =1}^\tau \Big[{1\over 1-{1\over p_\ell^{2\s}}} \Big]^{1/2}  ,\qquad \uz\in
\T^{\tau }.
\end{equation}

Thus
\begin{eqnarray}\label{uppersmalltau}\E\,  \sup_{\uz \in T^\tau } | \Upsilon(\uz)|&\le &    C_\s \prod_{\ell =1}^\tau \Big[{1\over 1-{1\over p_\ell^{2\s}}} \Big]^{1/2}  \Big({ \tau^{{3\over 2}- 2\s} \over(\log \tau)^{2 \s} }
    \Big)    . 
\end{eqnarray} 

\noi    Recall  that we have denoted $ \Pi_\s(\tau)= \prod_{\ell =1}^\tau  (  1-  p_\ell^{-2\s})^{-1} $.  By combining (\ref{uppersmalltau}) with (\ref{lowersmalltau}),  we get 
\begin{equation}\label{estsmalltau}
c_\s \, {\Pi_\s(\tau)^{1/2}\,\tau^{1-\s}  \over( \log \tau)^\s }\le \E \sup_{t\in \R}\big|\sum_{  n\le N\atop P^+(n)\le p_{\tau}}  {\e_n
\over n^{ \s+ it}}\big|
 \le C_\s \,    \Big({\Pi_\s(\tau)^{1/2}\, \tau^{{3\over 2}- 2\s} \over(\log \tau)^{2 \s} }
    \Big).
\end{equation}

 If $\s=1/2$, the modifications for $R$ and $S$ are, by using Mertens Theorem 
\begin{eqnarray}\label{debutR1demi} R   & \le &   C  \prod_{\ell =1}^\tau\Big[{1\over 1-{1\over p_\ell }} \Big]\Big( \sum_{j=1}^\tau {|z_{j } - z'_{j }|\over p_j }  
    \Big)^2  \le   C  (\log  \tau) \Big(\sum_{j=1}^\tau {1 \over p_j }  
    \Big) 
 \Big(\sum_{j=1}^\tau {|z_{j } - z'_{j }|^2\over p_j }  
    \Big) 
\cr & \le &  C  (\log  \tau)  (\log\log  \tau ) \Big(\sum_{j=1}^\tau {|z_{j } - z'_{j }|^2\over p_j }  
    \Big) 
 ,\end{eqnarray} 
and
  \begin{eqnarray}\label{debutS1demi}S 
  & \le &\sum_{j=1}^\tau|z_j - z'_j|
^2\sum_{b=1}^\infty {b^2\over p_j^{  b  }}\sum_{ m\le N/p_j^b   
  \atop P^+(m)\le p_\tau  } {1\over m  } 
  \le   C  \prod_{\ell =1}^\tau \Big[{1\over 1-{1\over p_\ell }} \Big]\sum_{j=1}^\tau|z_j - z'_j|
^2\sum_{b=1}^\infty {b^2\over p_j }
\cr& \le &  C  (\log  \tau) \sum_{j=1}^\tau   {|z_j - z'_j|
^2\over p_j  }  .
 \end{eqnarray}
Hence \begin{equation}  \label{estK1demi} \|\Upsilon(\uz)-\Upsilon(\uz)\big\|_2^2
 \le   C  (\log  \tau)  (\log\log  \tau ) \Big(\sum_{j=1}^\tau {|z_{j } - z'_{j }|^2\over p_j }  
    \Big) 
.
\end{equation}
And by Slepian's Lemma
$$\E\,  \sup_{\uz, \uz'\in T^\tau } |\Upsilon (\uz')-\Upsilon(\uz)|\le C  (\log  \tau)  (\log\log  \tau ) \Big(\sum_{j=1}^\tau {1\over p_j^{1/2} }  
    \Big)\le C\big({\tau  \log\log  \tau\over \log \tau }\big)^{1/2}.$$
 As 
 $$\|\Upsilon (\uz)\|_2 \le \big[ \sum_{n\in F_\tau  }{1\over n  }\big]^{1/2}\le C  \prod_{\ell =1}^\tau \big[{1\over 1-{1\over p_\ell }} \Big]^{1/2}\le C(\log \tau)^{1/2}  ,\qquad \uz\in
\T^{\tau },$$
we conclude to 
 \begin{equation} \label{e212demi}
\E\,  \sup_{\uz\in  T^\tau } |\Upsilon (\uz)|
\le   C    {\tau    }^{1/2} (\log\log  \tau )^{1/2} .\end{equation}
 Combining this estimate with (\ref{lowerdemismalltau}) finally gives
 \begin{equation} \label{e212demi.}
 C_1    {\tau    }^{1/2}  \le \E\,  \sup_{\uz\in  T^\tau } |\Upsilon (\uz)|
\le   C_2    {\tau    }^{1/2} (\log\log  \tau )^{1/2} .\end{equation}

\cqfd

\noi{\sl Acknowledgments:} {\it I thank Mikhail Lifshits for stimulating comments.}


{
}

 \noi {\phh Michel  Weber, \noi  Math\'ematique (IRMA),
Universit\'e Louis-Pasteur et C.N.R.S.,   7  rue Ren\'e Descartes,
67084 Strasbourg Cedex, France.
\par\noindent
E-mail: \  \tt weber@math.u-strasbg.fr} 
\end{document}